\documentclass[a4paper,11pt]{amsart}
\usepackage{amscd,amsthm,amsfonts,latexsym,amssymb}

\usepackage{diagrams}
\newarrow{Into} C--->
\newarrow{Mapsto} |--->
\newarrow{Corresponds} <--->

\newcommand{\ds}{\displaystyle}
\newtheorem{theorem}{Theorem} [section]
\newtheorem{lemma}[theorem]{Lemma}
\newtheorem{corollary}[theorem]{Corollary}

\newtheorem{proposition}[theorem]{Proposition}
\newtheorem{example}[theorem]{Example}

\newtheorem{definition}[theorem]{Definition}
\newtheorem{remark}[theorem]{Remark}

\newcommand{\inn}[1]{\langle#1\rangle} %inner product
\newcommand{\Biginn}[1]{\Big\langle#1\Big\rangle} %\Big version
\renewcommand{\vec}[1]{\boldsymbol{#1}} %for vectors (\vec usually
\newcommand{\dd}{{\mathrm d}}  %differential

\newcommand{\RR}{{\mathbb R}}  %real numbers
\newcommand{\CC}{{\mathbb C}}  %complex numbers

\newcommand{\MM}{{\mathbb M}}  %Minkowski
\newcommand{\DD}{{\mathbb D}}  %double (hyperbolic) numbers
\newcommand{\BB}{{\mathbb B}}  %bicomplex numbers

\newcommand{\Hh}{{\mathcal H}}  %horizontal space/projn
\newcommand{\Vv}{{\mathcal V}}  %vertical space or projn
\newcommand{\Dd}{{\mathcal D}}  %Degenerate subspaces 
\newcommand{\Kk}{{\mathcal K}}  %Complexified south pole
\newcommand{\Nn}{{\mathcal N}} %fattened origin

\newcommand{\Zz}{{\mathcal Z}} %null 1-d subspaces

\newcommand{\Q}{{\mathcal Q}}  %quadric
\newcommand{\CQ}{{\mathcal C}{\Q}} %cone on quadric

\newcommand{\RP}{{\mathbb R}P} %real projective space
\newcommand{\CP}{{\mathbb C}P} %complex projective space
\newcommand{\cc}{{\mathbf c}}
\newcommand{\ii}{{\rm i}}  %root of -1
\newcommand{\jj}{{\rm j}}  %hyperbolic number j

\newcommand{\pa}{\partial}

\renewcommand{\phi}{\varphi}
\newcommand{\la}{\lambda}
\newcommand{\La}{\Lambda}

\newcommand{\al}{\alpha}
\newcommand{\be}{\beta}
\newcommand{\ga}{\gamma}
\newcommand{\Ga}{\Gamma}
\newcommand{\ve}{\varepsilon}
\newcommand{\si}{\sigma}

\newcommand{\wt}{\widetilde} %usually better than \tilde
\newcommand{\ov}{\overline} %usually better than \ov
\newcommand{\SO}{\mbox{\rm SO}} %Special orthogonal group
\newcommand{\GL}{\mbox{\rm GL}} %General linear group
\newcommand{\C}{\mbox{\rm C}}   %Conformal group

\newcommand{\const}{\mbox{\rm const.}} %constant
\newcommand{\grad}{\mbox{\rm grad}\,}

\newcommand{\CN}{\mbox{\rm CN}}
\newcommand{\spn}{\mbox{\rm span}}

\begin{document}

\title{Harmonic morphisms and bicomplex manifolds}
\author{Paul Baird}
\author{John C.\ Wood*}
\thanks{The second author thanks the Gulbenkian foundation
 and the University of Lisbon for support and hospitality, and both authors thank the Agence National de Recherche, project ANR-07-BLAN-0251-01, for financial support.} 

\address{D\'epartement de Math\'ematiques, Universit\'e de Bretagne Occidentale, 6 Avenue Le Gorgeu, 29285 Brest, France}

\address{Department of Pure Mathematics, University of Leeds\\
Leeds LS2 9JT, Great Britain}
\email{Paul.Baird@univ-brest.fr; j.c.wood@leeds.ac.uk}

\begin{abstract} 
We use functions of a bicomplex variable to unify the existing constructions of harmonic morphisms from a $3$-dimensional Euclidean or pseudo-Euclidean space to a Riemannian or Lorentzian surface.  This is done by using the notion of complex-harmonic morphism between
complex-Riemannian manifolds and showing how these are given by bicomplex-holomorphic functions when the codomain is one-bicomplex dimensional.   By taking real slices, we recover well-known compactifications for the three possible real cases. On the way, we discuss some interesting conformal compactifications of complex-Riemannian manifolds by interpreting them as bicomplex manifolds.
\end{abstract}



\maketitle

\thispagestyle{empty}

\section{Introduction}   

\emph{Harmonic morphisms} are maps $\phi:M \to N$ between Riemannian or semi-Riemannian manifolds which preserve Laplace's equation in the sense that, if $f:V \to \RR$ is a harmonic function on an open subset of $N$ with $\phi^{-1}V$ non-empty, then $f \circ \phi$ is a harmonic function on $\phi^{-1}V$.
In the Riemannian case, they can be characterized as harmonic maps which are \emph{horizontally weakly conformal} (also called \emph{semiconformal}), a condition dual to weak conformality.  The characterization can be extended to harmonic morphisms between semi-Riemannian manifolds, with the additional feature that fibres can be degenerate.

Harmonic morphisms into Riemannian or Lorentzian surfaces are particularly nice as they are \emph{conformally invariant} in the sense that only the conformal equivalence class of the metric on the codomain matters; equivalently postcomposition of a harmonic morphism to a surface with a weakly conformal map of surfaces is again a harmonic morphism.
In particular, harmonic morphisms from (open subsets of) Minkowski $3$-space into $\CC$ are precisely the same as \emph{complex-valued null solutions of the wave equation}.

In \cite{Ba-Wo-Bernstein}, a twistorial Weierstrass-type representation was given which determined all harmonic morphisms from (convex) domains of $\RR^3$ to Riemann surfaces in terms of a pair of holomorphic functions; this led to a Bernstein-type theorem that the only globally defined harmonic morphism from $\RR^3$ to a Riemann surface is orthogonal projection, up to isometries and post-composition with weakly conformal maps.

In \cite{Ba-Wo-Rou}, a version of this was given for harmonic morphisms from Minkowski 3-space to Riemann surfaces, and also to \emph{Lorentz surfaces}, i.e., surfaces with a conformal equivalence class of metrics with signature $(1,1)$.  In the first case, the representation again involved holomorphic functions of a complex variable; however, in the second case, those were replaced by hyperbolic-holomorphic (`\emph{H-holo\-morphic}') functions of a variable which was a hyperbolic (also called `paracomplex') number $x+y\jj$ where $x$ and $y$ are real and $\jj^2 = 1$.  This led to interesting examples of globally defined harmonic morphisms other than orthogonal projection, and to harmonic morphisms all of whose fibres are degenerate.  In particular, it was shown that such degenerate harmonic morphisms correspond to
\emph{real-valued} null solutions of the wave equation.

\emph{Complex-Riemannian manifolds} were introduced by C.~LeBrun \cite{Le} as complex manifolds endowed with a \emph{symmetric} complex bilinear form on the holomorphic tangent space.  
The above constructions can be unified by employing \emph{complex-harmonic morphisms between complex-Riemannian manifolds}.  Complex-harmonic morphisms enjoy many of the properties of harmonic morphisms between semi-Riemannian manifolds, and have already been considered in \cite{Pa-Wo}, and by the authors in \cite{Ba-Wo-SFR}. 

The \emph{bicomplex numbers} are simultaneously a complexification of the complex numbers and the hyperbolic numbers. There is a natural notion of \emph{bicomplex-holomorphicity} which extends both holomorphicity and H-holo\-mor\-ph\-icity and leads to the notion of a \emph{bicomplex manifold}.

When the codomain of a complex-harmonic morphism has complex dimension 2, then we can consider it to be a
one-dimensional bicomplex manifold.  Our twistor data is bicomplex-holomorphic and naturally lives on such a manifold,
leading to our description in \S \ref{sec:cx-ha-bicx} of complex-harmonic morphisms from domains of $\CC^3$ into a one-dimensional bicomplex manifold. 

Once they have been given a `complex-orientation', the possible directions of non-degenerate fibres live in the complexification
$S^2_{\CC}$ of the $2$-sphere; to
allow degenerate fibres, we need the conformal compactification of this.  We describe that in two ways:
(i) the \emph{bicomplex quadric}
$\Q^1_{\BB}$\,,  (ii) the \emph{complex quadric} $\Q^2_{\CC}$ in $\CP^3$; we show that these are the same as one-dimensional bicomplex manifolds.
We obtain
an interesting diagram of compactifications, together with some double covers given by forgetting orientations, see
\S \ref{sec:interp}.  Finally, we show that all formulae and compactifications reduce to the known formulae and standard compactifications in the three real cases above.

One could extend this work to include harmonic morphisms from other $3$-dimensional space forms, treated in the Riemannian case in \cite{Ba-Wo-spaceform}, or to unify constructions of harmonic morphisms from suitable four-dimensional manifolds to surfaces, for example, Einstein anti-self-dual manifolds as in \cite{Wo-Einstein}.  This was partially done in \cite{Ba-Wo-SFR}  for Euclidean spaces by complexifying just the domain, showing that harmonic morphisms from $4$-dimensional Euclidean spaces to $\CC$ are equivalent to shear-free ray congruences or to Hermitian structures.

\section{Bicomplex numbers and bicomplex manifolds} \label{sec:bicomplex}

Bicomplex numbers have been invented and studied by many authors, often under a different name; a key paper is that of C.\ Segre in 1892 \cite{Se}.  The system of bicomplex numbers is the first non-trivial complex Clifford algebra
(and the only commutative one) and has recently been applied to quantum mechanics, see \cite{Ro-Tr1,Ro-Tr2} and the references therein, and to the study of Fatou and Julia sets in relation to $3$-dimensional fractals \cite{Ch-Ro-Sh} (see also the WEB page \cite{Rochon} for a list of related articles).  As will be seen below, bicomplex numbers give a natural way of complexifying formulae which are already complex, such as those which arise in twistor theory; we thus anticipate that more applications to physics will be found.  

In the sequel, we shall refer to \cite{Riley} and the more modern treatment given in \cite{Ro}.
The algebra of \emph{bicomplex numbers} is the space
$$
\BB = \{x_1+x_2\ii_1+x_3\ii_2+x_4\jj : x_1,x_2,x_3,x_4\in \RR\}. 
$$
As a real vector space, it is isomorphic to $\RR^4$ via the map 
\begin{equation} \label{B-R4}
\BB \ni x_1+x_2\ii_1+x_3\ii_2+x_4\jj \mapsto (x_1,x_2,x_3,x_4) \in \RR^4,
\end{equation}
from which it inherits its additive structure.  Multiplication is defined by the rules:
$$
\ii_1{}\!^2=\ii_2{}\!^2 = - 1, \qquad \ii_1\ii_2=\ii_2\ii_1=\jj \quad \text{so that} \quad \jj^2 = 1\,.
$$
Let $\CC [\ii_1]$ denote the field of complex numbers $\{ x + y \ii_1: x,y\in \RR\}$. We can write any $q \in \BB$ as
\begin{equation} \label{q-z1-z2}
q = q_1 + q_2 \ii_2 \quad \text{where } q_1,
	q_2\in \CC[\ii_1]\,;
\end{equation}
comparing with \eqref{B-R4} we have $q_1 = x_1+x_2 \ii_1$ and $q_2 = x_3+ x_4 \ii_1$; 
the map $q \to (q_1,q_2)$ gives a natural isomorphism between the vector spaces $\BB$ and $\CC^2$.
With the notation \eqref{q-z1-z2}, multiplication takes the form
$$
(q_1+q_2\ii_2)(w_1+w_2\ii_2) = q_1w_1-q_2w_2 + (q_1w_2+q_2w_1)\ii_2\,;
$$
thus $\BB$ can be viewed as a natural extension of the complex number system $\CC[\ii_2] = \{ x + y \ii_2 : x,y\in \RR\}$, but now with $x,y \in \CC[\ii_1]$; in other words, $\BB = \CC \otimes_{\RR} \CC$.  However, unlike the complex numbers, the algebra $\BB$ has
\emph{zero divisors}, namely the set of points 
$\{ q_1+q_2\ii_2 \in \BB : q_1{}\!^2 + q_2{}\!^2 = 0\}$ $= \{ z(1 \pm \jj) : z \in \CC[\ii_1]\}$\,.
Following \cite{Ro}, we call the complex number $\CN(q):= q_1{}\!^2 + q_2{}\!^2$ the \emph{complex (square) norm of} $q$.
Then a bicomplex number $q = q_1+q_2\ii_2$ is a \emph{unit}, i.e., has an inverse, if and only if $\CN(q) \neq 0$; its inverse is then given by
$q^{-1} = (q_1 - q_2 \ii_2) \big/\CN(q)$\,. The set of units forms a multiplicative group which we denote by $\BB_*$\,.
Writing $q^* = q_1 - q_2\ii_2$, we see that $\CN(q) = qq^*$;  hence, if $\CN(q) \neq 0$, then $q^{-1} = q^*/\CN(q)$.
Note that $q$ also inherits a \emph{real norm} from $\RR^4$ given by
$|q| = \sqrt{|q_1|^2 + |q_2|^2}
 = \sqrt{x_1^{}\!^2 + x_2^{}\!^2 + x_3^{}\!^2 + x_4^{}\!^2}$\,.
 
By the \emph{positive complex conformal group} we mean the matrix group  
 \begin{equation} \label{C+}
\C_+(2,\CC) = \{A \in \GL(2,\CC) : A^T A = (\det A)I \}\,;
\end{equation}
this is the identity component of the complex conformal group $\{A \in \GL(2,\CC) : A^T A = \lambda I \text{ for some } \lambda \in \CC \}$, noting that $\lambda$ is necessarily $\pm\det A$.  
The map
$q =q_1+q_2\ii_2 \mapsto \left(\begin{smallmatrix}q_1 & -q_2 \\ q_2 & q_1 \end{smallmatrix}\right)$
is an algebra-homomorphism from $\BB$ to the $2 \times 2$  complex matrices with the group $\BB_*$ of units mapping onto
$\C_+(2,\CC)$, and the bicomplex numbers of complex norm one mapping onto the complex special orthogonal group
$\SO(2,\CC) = \{ A \in \GL(2,\CC) : \det A = 1, A^T A =I\}$.

We generalize these notions to \emph{bicomplex vectors}
$\vec{q} = (q_1, \ldots, q_m)\in \BB^m$.
Extend the standard com\-plex-bi\-lin\-ear inner product $\inn{ \ , \ }_{\CC}$
to a bicomplex-bilinear inner product $\inn{ \ , \ }_{\BB}$ on
$\BB^m$; explicitly, for
$\vec{p} = (p_1, \ldots, p_m) \in \BB^m$, we have
$\inn{\vec{p}, \vec{q}}_{\BB} = \sum_{k=1}^m p_k q_k$.  Then, for
a bicomplex vector $\vec{q} = \vec{u} + \vec{v}\ii_2$ \ $(\vec{u}, \vec{v} \in \CC[\ii_1]^m)$ we have four important quantities:

(i) the bicomplex number $\vec{q}^2 := \inn{\vec{q},\vec{q}}_{\BB} = \sum_{k=1}^m q_k{}\!^2$.
Note that
$\vec{q}^2  =
\inn{\vec{u}+\vec{v}\ii_2, \vec{u}+\vec{v}\ii_2}_{\BB} = \inn{\vec{u}, \vec{u}}_{\CC} - \inn{\vec{v}, \vec{v}}_{\CC} + 	2\inn{\vec{u},\vec{v}}_{\CC}\ii_2$\,;

(ii) the bicomplex vector $\vec{q}^* = (q_1{}\!^*, \ldots , q_m{}\!^*)$;

(iii) the \emph{complex (square) norm}  $\CN(\vec{q}):= \vec{q}\vec{q}^* = \sum_{k=1}^m\CN(q_k) \in \CC$. We have
$\CN(\vec{q}) = \vec{u}^2 + \vec{v}^2$ where we write $\vec{u}^2 = \inn{\vec{u},\vec{u}}_{\CC}$ and $\vec{v}^2 = \inn{\vec{v},\vec{v}}_{\CC}$. Note that
 $\CN(\la \vec{q}) = \CN(\la ) \CN(\vec{q})$ for $\la \in \BB$;

(iv) \emph{the real norm}
$|\vec{q}| = \sqrt{\sum_{k=1}^m |q_k|^2}
= \sqrt{|\vec{u}|^2 + |\vec{v}|^2}$,
which we only use for notions of convergence.

The complex numbers embed naturally in $\BB$ via the inclusion:
\begin{equation} \label{embed-cx}
\iota_{\CC}: \CC \hookrightarrow \BB\,,  \quad \iota_{\CC} (x+y\ii ) = x + y \ii_2  \quad (x,y \in \RR)\,;
\end{equation}
the use of $i_2$ rather than $i_1$ is a convention which carries through to all our formulae.
However, the alternative embedding $z = x+ y \ii \mapsto x+y \ii_1 = z + 0 \ii_2$ appears in various places including Example \ref{ex:BaWo2}. 

Now let $ \phi : U \to \BB$ be a function defined on an open subset of $\BB$; write 
\begin{equation} 
 \psi(q_1+q_2\ii_2) =  \psi_1(q_1,q_2) +  \psi_2(q_1, q_2) \ii_2\,.
\end{equation}
Here we take $\psi_1$ and $\psi_2$ to be holomorphic in $(q_1, q_2)$ --- this turns out to be a necessary condition for the existence of the bicomplex derivative which we now define. Specifically, let $p \in U$.  Then the \emph{bicomplex derivative of the function $q \mapsto \psi(q)$ at $p$} is the limit
$$
\psi^{\prime}(p) : = \frac{\dd\psi}{\dd q}(p)
: = \lim_{|h|\to 0,\; \CN(h) \neq 0} \frac{\psi(p+h) - \psi (p)}{h}\;,
$$
whenever this exists.  It is easy to see that the bicomplex derivative of $\psi = \psi_1 + \psi_2 \ii_2$ exists if and only if the pair $(\psi_1,\psi_2)$ of holomorphic functions satisfies the following \emph{bicomplex Cauchy-Riemann equations}:
$$
\dfrac{\pa \psi_1}{\pa q_1}  =  \dfrac{\pa \psi_2}{\pa q_2}
\quad \text{and} \quad
\dfrac{\pa \psi_1}{\pa q_2}  =  - \dfrac{\pa \psi_2}{\pa q_1}\,.
$$
When this is the case, we shall say that $\psi$ is
\emph{bi\-com\-plex-differ\-enti\-able} or \emph{bi\-com\-plex-holo\-mor\-phic}.  Note that the bicomplex Cauchy-Riemann equations are equivalent to the condition that the differential of $\psi$ lie in $\C_+(2,\CC)$.

On defining partial derivatives formally by
$$
\frac{\pa\psi}{\pa q}
= \frac{1}{2}\Bigl(\frac{\pa\psi}{\pa q_1}
	- \frac{\pa\psi}{\pa q_2}\ii_2 \Bigr)\,, \quad 
\frac{\pa\psi}{\pa q^*}
= \frac{1}{2}\Bigl(\frac{\pa\psi}{\pa q_1}
	+ \frac{\pa\psi}{\pa q_2}\ii_2 \Bigr)
$$
where $\pa\psi/\pa q_k
= \pa\psi_1/\pa q_k + (\pa\psi_2/\pa q_k)\ii_2$ \ $(k = 1,2)$,
the bicomplex Cauchy-Riemann equations can be written as the single equation: $\pa\psi/\pa q^* = 0$.

Under the embedding \eqref{embed-cx}, holomorphic maps extend to bi\-com\-plex-holo\-mor\-phic maps as follows, the proof is by analytic continuation.

\begin{lemma} \label{lem:derivs}
Let $f:U \to \CC$ be holomorphic map from an open subset of\/ $\CC$. 
Then $f$ can be extended to a bi\-com\-plex-holo\-mor\-phic function $\psi:\wt{U} \to \BB$ on an open subset $\wt{U}$ of\/ $\BB$ with $\wt{U} \cap \CC = U$; the germ of the extension at $U$ is unique. 

Conversely, the restriction of any
bi\-com\-plex-holo\-mor\-phic function $\wt{U} \to \BB$ to
$U = \wt{U} \cap \CC$ is holomorphic,
provided that $U$ is non-empty.
\qed \end{lemma} 

\begin{remark} \label{rem:Ringleb} \rm
Another way to understand bi\-com\-plex-holo\-mor\-phic functions is \emph{Ringleb's Lemma} \cite[\S 9]{Riley} as follows.  Let $a = \frac{1}{2}(1-\jj)$ and $b = \frac{1}{2}(1+\jj)$; then $a$ and $b$ are zero divisors with $a^2=a$, $b^2=b$ and $ab=0$. 

Any bicomplex number $q \in \BB$ can be written uniquely in the form $q = z a + w b$ with $z, w \in \CC[\ii_1]$
thus identifying $\BB$ with $\CC\oplus\CC$.  Then $\psi$ is bi\-com\-plex-holo\-mor\-phic if and only if it is of the form $\psi(q) = f_1(z) a + f_2(w) b$ for some holomorphic functions $f_1$ and $f_2$. 

With this formulation, a biholomorphic function $\psi$ is an extension of a holomorphic function $f:U \to \CC$ if and only if $f_1 = f_2 = f$.
\end{remark}

By a \emph{bicomplex manifold} we mean a complex manifold with a complex atlas whose transition functions are bicomplex-holomorphic functions.  Such a complex manifold is necessarily of even dimension $2n$; we call $n$ the
\emph{bicomplex dimension}.  Then a map between bicomplex manifolds is called
\emph{bicomplex-holomorphic} if it is bicomplex-holomorphic in all the charts.  Note that, by Ringleb's lemma, a bicomplex manifold of bicomplex dimension $n$ is locally the product of complex manifolds of dimension $n$; however, the complex-Riemannian metrics we introduce later are never product metrics, so this observation is of limited use in our work.

Bicomplex manifolds can be obtained by complexifying complex manifolds; we give some examples that we shall use later.

\begin{example} \label{ex:S2C} {\rm (Complex $2$-sphere)}
The \emph{complex $2$-sphere} is the complex surface
$$
S^2_{\CC} = \{(z_1,z_2,z_3) \in \CC^3: z_1^{}\!^2 + z_2^{}\!^2 +z_3^{}\!^2 = 1\};
$$
this may be considered as a complexification of the usual $2$-sphere $S^2$.
We give some charts.

{\rm (i)} Set $\Hh^1 = \{G \in \BB : \CN(G) = -1\}$, and $\Kk^1 = \{(z_1,z_2,z_3) \in S^2_{\CC}: z_1 = -1\}$, the `complexified' south pole.
We have a bijection $\si_{\CC}:U_G \to \BB \setminus \Hh^1$, $(z_1,z_2,z_3) \mapsto (z_2 +  z_3\ii_2)/(1+z_1)$ from
$U_G = S^2_{\CC} \setminus \Kk^1$, with inverse
\begin{equation} \label{stereo}
G = G_1 + G_2\ii_2  \mapsto
\bigl(1-\CN(G)\,,\, 2G_1\,,\,2G_2 \bigr) \big/ \bigl(1+\CN(G) \bigr)\,;
\end{equation}
note that this is the complexification of standard stereographic projection on $S^2 \setminus \{(0,0,-1)\}$.
We call this the \emph{standard chart} for the complex $2$-sphere.

$\check{({\mathrm i})}$
Similarly, stereographic projection from the north pole complexifies to give a bijection
$\check{\si}_{\CC}: U_{\check{G}} \to \BB \setminus \Hh^1$ where $U_{\check{G}} = S^2_{\CC} \setminus \check{\Kk}^1$ with
$\check{\Kk}^1 = \{(z_1,z_2,z_3) \in S^2_{\CC}: z_1 = +1\};$
this has inverse
$$
\check{G} \mapsto \bigl(\CN(\check{G}) - 1\,,\, 2\check{G}_1\,,\, -2\check{G}_2 \bigr) \big/ \bigl( \CN(\check{G}) + 1 \bigr).
$$
These two charts cover $S^2_{\CC}$, i.e.,
$U_G \cup U_{\check{G}} = S^2_{\CC}$.
Further,
$\si_{\CC}(U_G \cap U_{\check{G}})
= \check{\si}_{\CC} (U_G \cap U_{\check{G}})
= \BB_* \setminus \Hh^1$  and the transition function
$\check{\si}_{\CC} \circ \si_{\CC}^{-1}:\BB_* \setminus \Hh^1 \to
\BB_* \setminus \Hh^1$ is
$\check{G} = 1/G,$ so that the two charts give $S^2_{\CC}$ the structure of a one-dimensional bicomplex manifold.

Many other bicomplex charts can be obtained  by simple modifications of these; for comparison with other spaces we shall need the following$:$

{\rm (ii)} $L = L_1 + L_2 \ii_2 \mapsto (-2L_2,1- \CN(L), -2L_1) \big/ (1+\CN(L))$
defines a chart which maps  $\BB \setminus \Hh^1$  to
$S^2_{\CC} \setminus \{(z_1,z_2,z_3) \in S^2_{\CC} : z_2 = -1\};$

{\rm (iii)} $K = K_1 + K_2 \ii_2 \mapsto (-2K_1, -2K_2, 1-\CN(K)) \big/ (1+\CN(K))$
defines a chart which maps  $\BB \setminus \Hh^1$  to
$S^2_{\CC} \setminus \{(z_1,z_2,z_3) \in S^2_{\CC} : z_3 = -1\}$.

The transition functions with the standard chart are
\begin{align} 
L &= (G-1)\ii_2\big/ (G+1)   &\text{with  inverse} &
	&G = (1 - L\ii_2 ) \big/ (1+ L\ii_2),
\label{transitions-L}\\
K &= (G - \ii_2)/(G+\ii_2)   &\text{with inverse}  &
	&G= (1+K)\ii_2/(1 - K) \label{transitions-K}.
\end{align}
Both of these maps are bicomplex-holomorphic functions on their domains with
bicomplex-holomorphic inverses; their domains and ranges are easily calculated, for example, 
\eqref{transitions-L} is a bijection from
\newline
$\BB \setminus \Hh^1 \setminus \{G \in \BB: \CN(1+G) = 0\}$  to
$\BB \setminus \Hh^1 \setminus \{K \in \BB: \CN(1-K) = 0\}$.
\end{example}

The next two examples are less obvious.

\begin{example} \label{ex:bicx-quadric} {\rm (Bicomplex quadric)}
Let $\Nn$ be the `fattened origin'\\
$\Nn = \{\vec{\xi} \in \BB^3 : \CN(\xi_i) = 0 \ \forall i=1,2,3 \}$, and let
$$
\CQ^1_{\BB}
= \{\vec{\xi} \in \BB^3 \setminus \Nn : \vec{\xi}^2 = 0\}
	= \bigl\{(\xi_1, \xi_2, \xi_3) \in \BB^3 \setminus \Nn :
\xi_1^{}\!^2 + \xi_2^{}\!^2 + \xi_3^{}\!^2 = 0\bigr\}.
$$
Define an equivalence relation on $\CQ^1_{\BB}$ by
$\vec{\xi} \sim \vec{\tilde{\xi}}$ if $\vec{\tilde{\xi}} = \la \vec{\xi}$
for some $\la \in \BB;$ note that $\la$ is necessarily a unit, for otherwise
 $\vec{\tilde{\xi}}$ would lie in $\Nn$.  We call the set of equivalence classes
the \emph{bicomplex quadric} $\Q^1_{\BB}$.
We can give this the structure of a one-dimensional bicomplex manifold by using the following charts  which cover $\Q_{\BB}^1$\,.
 
{\rm (i)} $G \mapsto [-2G, 1- G^2, (1+G^2)\ii_2]$ maps $B_*$ onto the open set\\
$U_G = \{ [\vec{\xi}]\in \Q^1_{\BB}: \CN(\xi_1) \neq 0\} $ and has inverse
\begin{equation} \label{G-inv}
G  = (\xi_2 + \xi_3\ii_2) \big/ \xi_1 = -\xi_1 \big/ (\xi_2 - \xi_3\ii_2) \,.
\end{equation}
Note that $\CN(\xi_1) \neq 0$ implies that
$\CN(\xi_2 - \xi_3 \ii_2) \neq 0$  and $\CN(\xi_2 +\xi_3 \ii_2) \neq 0$
from the following \emph{fundamental identity} valid for all
$\vec{\xi} \in \BB^3$ with $\vec{\xi}^2 = 0:$
$$
\CN(\xi_1)^2 =  \CN(\xi_2 - \xi_3 \ii_2)\,\CN(\xi_2 + \xi_3 \ii_2)\,;
$$
thus both fractions in \eqref{G-inv} are well-defined; we see easily that $\CN(G) \neq 0$.
We shall refer to this chart as the \emph{standard chart}.

$\check{({\mathrm i})}$  The chart $\check{G} \mapsto [-2\check{G},  \check{G}^2-1, (\check{G}^2+1)\ii_2]$ maps $B_*$ onto the \emph{same} open set $U_G$.  Note that the transition function with the standard chart is $\check{G} = 1/G$ on $B_*$, as before.

{\rm (ii)} $L \mapsto [(1+L^2)\ii_2, 2L, 1-L^2]$ maps $B_*$ onto the open set $U_{L} = \{ [\vec{\xi}]\in \Q^1_{\BB} : \CN(\xi_2) \neq 0\}$ and has inverse
$
L = -(\xi_3 + \xi_1\ii_2) \big/ \xi_2 = \xi_2 \big/ (\xi_3 - \xi_1\ii_2)\,.
$

{\rm (iii)} $K \mapsto [1 - K^2, (1+K^2)\ii_2, 2K]$ maps $B_*$ onto the open set  $U_K  = \{ [\vec{\xi}]\in \Q^1_{\BB} : \CN(\xi_3) \neq 0\}$ and has inverse
$
K = -(\xi_1 + \xi_2\ii_2) \big/ \xi_3 = \xi_3 \big/ (\xi_1 - \xi_2\ii_2)\,.
$

Clearly $U_G\cup U_{L} \cup U_K = \Q^1_{\BB}$.  It can be checked that the transition functions  are given by \eqref{transitions-L}
and \eqref{transitions-K}  on suitable domains.
Since these are  bicomplex-holomorphic, the three charts give the bicomplex quadric the structure of a one-dimensional bicomplex manifold.

Let $\Q^1_{\CC} = \bigl\{ \vec{\xi} = [ \xi_0,\xi_1,\xi_2] \in \CP^2 : \xi_0^{}\!^2 + \xi_1^{}\!^2 + \xi_2^{}\!^2 = 0 \bigr\} \cong \CP^1$.  With $a$ and $b$ as in Remark \ref{rem:Ringleb}, it can be checked that the map $\bigl([\vec{\eta}] , [\vec{\rho}] \bigr) \mapsto [\vec{\eta} a + \vec{\rho} b] $ is well-defined and gives a
complex diffeomorphism from $\Q^1_{\CC} \times \Q^1_{\CC}$ to $\Q^1_{\BB};$ thus as a complex manifold,
$Q^1_{\BB}$ is the product\/ $\CP^1 \times \CP^1$.
Note also that $\Q^1_{\BB}$ has a dense open subset
\begin{equation} \label{Q1B*}
\Q^1_{\BB *}  = \{[\vec{\xi}] \in \Q^1_{\BB} : \CN(\vec{\xi}) \neq 0\}\,;
\end{equation}
this is a one-dimensional bicomplex manifold which is not globally a product of complex curves.
\end{example}

\begin{example} \label{ex:complex-quadric} {\rm (Complex quadric)}
Let
$$
\Q^2_{\CC} = \bigl\{ [\zeta_0,\zeta_1,\zeta_2,\zeta_3] \in \CP^3: \zeta_0{}\!^2 = \zeta_1{}\!^2 + \zeta_2{}\!^2 + \zeta_3{}\!^2\bigr\};
$$ 
the choice of signs is the most convenient for later comparison with real cases, but is unimportant here.
This is again a one-dimensional bicomplex manifold.
Indeed the following maps give charts which cover $\Q^2_{\CC}$;
in formulae {\rm (i)} and ${\rm \check{(i)}}$, for convenience of notation, we identify the last two components $(\zeta_2, \zeta_3)$ of points of $\Q^2_{\CC}$ with the bicomplex number $\zeta_2 + \zeta_3\ii_2$.

{\rm (i)} $G \mapsto [1+\CN(G), 1-\CN(G), 2G]$ maps $\BB$ onto the 
open set $V_G = \{ [\vec{\zeta}] \in \Q^2_{\CC} : \zeta_0 + \zeta_1 \neq 0\} $
and has  inverse
$$
G = (\zeta_2 + \zeta_3\ii_2) \big/ (\zeta_0 + \zeta_1) \,.
$$
We shall refer to this as the \emph{standard chart} for $\Q^2_{\CC}$.   

$\check{({\mathrm i})}$ $\check{G} \mapsto [1+\CN(\check{G}), \CN(\check{G}) - 1, 2\check{G}^*]$ maps $\BB$ onto the open set
$V_{\check{G}} = \{ [\vec{\zeta}]\in \Q^2_{\CC}: \zeta_0 - \zeta_1 \neq 0\}$ and has inverse
$$
\check{G} = (\zeta_2 - \zeta_3 \ii_2) \big/ (\zeta_0 - \zeta_1) \,.
$$
The transition function with the standard chart is again $G = 1/\check{G}$ on $B_*$.

Both of these charts miss out the points $[0, 0, 1, \pm \ii_1]\in \Q^2_{\CC}$ so
we require another chart.  This can be either of the following.

{\rm (ii)} $L = L_1 + L_2 \ii_2 \mapsto [1+\CN(L), - 2L_2, 1 - \CN(L), - 2L_1]$ maps $\BB$ to the open set 
$V_L = \{ [\vec{\zeta}] \in \Q^2_{\CC} : \zeta_0 + \zeta_2 \neq 0\}$
and has inverse
$$
L = -(\zeta_3 + \zeta_1 \ii_2) \big/ (\zeta_0 + \zeta_2).
$$

{\rm (iii)}  $K = K_1 + K_2 \ii_2 \mapsto [1 + \CN(K), - 2K_1, - 2K_2, 1 - \CN(K)]$  maps $\BB$ to the  open set
$V_K = \{ [\vec{\zeta}] \in \Q^2_{\CC} : \zeta_0 + \zeta_3 \neq 0\}$
and has inverse
$$
K = -(\zeta_1 + \zeta_2 \ii_2) \big/ (\zeta_0 + \zeta_3).
$$

Again it can be checked that the transition functions are given by
\eqref{transitions-L} and \eqref{transitions-K} on suitable domains.

Note that $S^2_{\CC}$ embeds into $\Q^2_{\CC}$ via the mapping
$(z_1, z_2, z_3) \mapsto [1, z_1, z_2, z_3];$ this is clearly bicomplex-holomorphic.
\end{example}

We shall see later that the bicomplex quadric $\Q^1_{\BB}$ and the complex quadric $\Q^2_{\CC}$ are, in fact, \emph{equivalent} as bicomplex manifolds.  However, $S^2_{\CC}$ is \emph{not} compact, but has conformal compactification given by $\Q^1_{\BB} \cong \Q^2_{\CC}$\,.  

\smallskip

\section{Complex-harmonic morphisms} \label{sec:complex}

Let $U$ be an open subset of $\CC^m$.
Then we say that a holomorphic function $f : U \to \CC$ is \emph{com\-plex-har\-mon\-ic} if it satisfies the
\emph{complex-Laplace equation}:
$$
\Delta_{\CC}f : = \sum_{k = 1}^m\frac{\pa^2 f}{\pa z_i{}\!^2} = 0\,,
$$
where $(z_1, \ldots , z_m)$ are the
standard coordinates on $\CC^m$. 

More generally, let $M$ be a complex manifold of some complex dimension $m$; denote its
$(1,0)$- (holomorphic) tangent space by $T'M$; thus $T'M$ is spanned by $\{\pa/\pa z^i : i=1,\ldots, m\}$ for any complex coordinates $(z^i)$.  Following C.~LeBrun \cite{Le}, a holomorphic section $g$ of $T'M \otimes T'M$ which is symmetric and non-degenerate is called a \emph{holomorphic metric}; the pair $(M,g)$ is then called a
\emph{com\-plex-Riemann\-ian manifold}. The simplest example is the complex manifold $\CC^m$ endowed with its standard holomorphic metric $g = \dd z_1{}\!^2 + \cdots + \dd z_m{}\!^2$; this can be thought of as the complexification of $\RR^m$ with its standard metric.
More generally, if $(M_{\RR},g_{\RR})$ is a real-analytic Riemannian or semi-Riemannian manifold, then it has a germ-unique complexification $M_{\CC}$ with holomorphic tangent bundle $T'M_{\CC} = TM_{\RR} \otimes_{\RR} \CC$; extending the Riemannian metric by complex bilinearity to $T'M_{\CC}$ gives a holomorphic metric. 
For example, complexifying the $2$-sphere $S^2$ with its standard Riemannian metric gives the complex-Riemannian manifold $(S^2_{\CC}, g)$ with $g$ equal to the restriction of the standard holomorphic metric on $\CC^3$.  

A holomorphic function $f : M \to \CC$ from a complex-Riemannian manifold is said to be
\emph{com\-plex-har\-mon\-ic} if it satisfies the \emph{complex-Laplace equation} $\Delta^M_{\CC}f = 0$
where the complex-Laplace operator $\Delta^M_{\CC}$ is defined by complexifying the  formulae for the real case, for example, in local complex coordinates $(z^i)$, defining the matrix $(g_{ij})$ by
 $g_{ij} = g(\pa/\pa z_i, \pa/\pa z_j)$ and letting $(g^{ij})$ denote its inverse, we have
$$
\Delta^M_{\CC}f = g^{ij} \left(\frac{\pa^2 f}{\pa z^i
 \pa z^j} - \Ga_{ij}^k \frac{\pa f}{\pa z^k} \right)
\text{ where } \Ga_{ij}^k
	= \frac{1}{2}g^{km}\Bigl\{\frac{\pa g_{jm}}{\pa z_i} 
	+ \frac{\pa g_{im} }{\pa z_j}
		- \frac{\pa g_{ij}}{\pa z_m} \Bigr\}.
$$

\begin{definition} Let $(M,g)$ and $(N,h)$ be complex-Riemannian manifolds.  A holomorphic mapping $\Phi : M \to N$ is a \emph{com\-plex-har\-mon\-ic morphism} if, for every com\-plex-har\-mon\-ic function $f : V \to \CC$ defined on an open subset $V$ of\/ $N$ such that $\Phi^{-1}(V)$ is non-empty, the composition $f\circ \Phi : \Phi^{-1}(V) \to \CC$ is com\-plex-har\-mon\-ic.
\end{definition}

Clearly, many notions and results for harmonic morphisms between semi-Riemannian manifolds complexify immediately to complex-harmonic morphisms between complex-Riemannian manifolds. 
In particular, given a holomorphic map $\phi : (M, g) \to (N,h)$ between complex-Riemannian manifolds, its differential $\dd\phi_p:T_p'M \to T_{\phi(p)}'N$ at a point $p \in M$ is a complex linear map between holomorphic tangent spaces.  We
say that a holomorphic map 
$\phi:(M^m,g) \to (N^n,h)$ is \emph{complex-weakly conformal  with (complex-) square conformality factor} $\La (p)$ if
\begin{equation} \label{cx-WC}
h(\dd\phi_p(X),\dd\phi_p(Y)) = \Lambda(p)\,g(X,Y) \qquad (p \in M^m, \ X,Y \in T_p  M^m)
\end{equation}
for some holomorphic function $\Lambda:M^m \to \CC$.
In local complex coordinates, this reads
$$
h_{\al\be} \frac{\pa\phi^{\al}}{\pa z^i} \frac{\pa\phi^{\be}}{\pa z^j} = \Lambda\, g_{ij}\,.
$$

However, it is the following dual notion which is  more important to us.  We call $\phi$ \emph{(complex-) horizontally (weakly) conformal (complex-HWC)} with
\emph{(complex-)square dilation} $\La (p)$ if
\begin{equation} \label{HWC}
g\bigl(\dd \phi_p^*(U), \dd \phi_p^*(V)\bigr) = \La (p)\, h(U,V) \qquad (p \in M^m, \ U,V\in T_{\phi(p)}'N)
\end{equation}
for some holomorphic function $\Lambda:M^m \to \CC$ where $\dd \phi_p^*: T_{\phi(p)}'N \to T_p'M$ denotes the adjoint of $\dd\phi_p$ with respect to $g$ and $h$.  In local complex coordinates this reads
$$
g^{ij} \frac{\pa\phi^{\al}}{\pa z^i}\frac{\pa\phi^{\be}}{\pa z^j}
	= \Lambda h^{\al\be}.
$$

A subspace $W$ of $T_p'M$ is called \emph{degenerate} if there exists a non-zero vector $v\in W$ such that $g(v,w) = 0$ for all $w\in W$, and \emph{null} if $g(v,w) = 0$ \emph{for all} $v, w\in W$.  As in the semi-Riemannian case (see \cite[Proposition 14.5.4]{book}), a complex-HWC map can have three types of points, as follows; we use
${}^{\perp_{\cc}}$ to denote the orthogonal complement of a subspace in $T'M$ with respect to $g$.

\begin{proposition} \label{prop:types}
Let $\phi:(M,g) \to (N,h)$ be a complex-HWC map.  Then, for each $p \in M$, precisely one of the following holds$:$
 
{\rm (i)} $\dd\phi_p = 0$.  Then $\Lambda(p) = 0;$

{\rm (ii)} $\La(p) \neq 0$.  Then $\phi$ is submersive at $p$ and $\dd\phi_p$ maps the \emph{complex-horizontal space}
 $\Hh^{\cc}_p := (\ker\dd\phi_p)^{\perp_{\CC}}$ conformally onto $T_{\phi(p)}'N$ with square conformality factor $\La(p)$, i.e.,
$h(\dd\phi_p(X), \dd\phi_p(Y)) = \La(p)\,g(X,Y)$ \ $(X,Y \in \Hh_p),$
we call $p$ a \emph{regular point} of\/ $\phi;$

{\rm (iii)} $\La(p) = 0$ but $\dd\phi_p \neq 0$.  Then the \emph{vertical space} $\Vv_p^{\cc} := \ker\dd\phi_p$ is degenerate and $\Hh_p^{\cc} \subseteq \Vv_p^{\cc};$ equivalently, $\Hh_p$ is null and non-zero.  We say that $p$ is a \emph{degenerate point} of\/ $\phi$, or that $\phi$ is \emph{degenerate at $p$}.
\qed \end{proposition}

Other useful results are: (i) \emph{if $M$ and $N$ are complex surfaces}, by which we mean complex-Riemannian manifolds of complex dimension $2$,
\emph{a holomorphic map $\phi:M \to N$ is a harmonic morphism if and only if it is complex-HWC}.
As in the semi-Riemannian case,
see \cite[Remark 14.5.7]{book}, this condition is \emph{not} equivalent to complex-weakly conformal ---  behaviour at degenerate points is different;
(ii) \emph{the composition of a complex-harmonic morphism to a complex surface with a complex-HWC map of complex surfaces is another complex-harmonic morphism};
(iii) \emph{the concept of complex-harmonic morphism to a complex surface depends only on the conformal class of its holomorphic metric}.

We extend the fundamental characterization of harmonic morphisms between Riemannian or semi-Riemannian manifolds as horizontally weakly conformal harmonic maps \cite{Fu,Fu-2,Is} to the case of interest to us.
We use the standard com\-plex-bi\-lin\-ear inner product $\inn{ \ , \ }_{\CC}$ on $\CC^m$ and the complex gradient $\grad_{\CC}f =
\bigl( \pa f \big/ \pa z_1\,,\ldots,\, \pa f \big/ \pa z_m \bigr)$ of a holomorphic function $f$ defined on a subset of $\CC^m$.

\begin{proposition} \label{prop:fund-char} {\rm (Fundamental characterization)} Let $(M^m,g)$ be a complex-Riemannian manifold.
A holomorphic map $\Phi : M^m \to \CC^n$ is a com\-plex-har\-mon\-ic morphism if and only if it is complex-harmonic and complex-HWC;
explicitly, on writing $\Phi = ( \Phi_1, \ldots , \Phi_n)$, we have
\begin{equation} \label{cx-char}
\left\{
\begin{array}{lll} {\rm (a) } & \Delta_{\CC} \Phi_{\al} = 0 & (\al = 1, \ldots , n)\,, \\
 {\rm (b) } & \inn{\grad_{\CC}\Phi_{\al}, \grad_{\CC}\Phi_{\be}}_{\CC} = \delta_{\al\be}\La & (\al, \be = 1, \ldots , n)\,,
\end{array} \right.
\end{equation}
for some (holomorphic) function $\La : M^m \to \CC$. 
\end{proposition}

\begin{proof}  Suppose that $\Phi$ is a com\-plex-har\-mon\-ic morphism. Given a point $p \in \CC^n$ and complex constants
$\{ C_{\al}, C_{\al\be}\}_{\al, \be = 1, \ldots , n}$ with $\C_{\al \be} = C_{\be \al}$ and $\sum_{\al = 1}^nC_{\al\al} = 0$, then, writing $(w_1, \ldots , w_n)$ for the standard complex coordinates on $\CC^m$, there exists a com\-plex-har\-mon\-ic function $f$ defined on a neighbourhood of $p$ with 
$$
\frac{\pa f}{\pa w_{\al}}(p) = C_{\al}  \quad {\rm and} \quad \frac{\pa^2f}{\pa w_{\al}\pa w_{\be}} (p) = C_{\al \be} \quad (\al , \be = 1, \ldots , n)\,;
$$
we simply take $f = C_{\al \be} w_{\al}w_{\be} + C_{\al}w_{\al}$
(summing over repeated indices).

Now, let $p \in M^m$ and let $z^i$ be local complex coordinates on a neighbourhood of $p$ such that $g^{ij} = \delta_{ij}$
at $p$.  Then, by the composition law,
\begin{equation} \label{chain3}
\Delta_{\CC} (f \circ \Phi )
= \frac{\pa f}{\pa w_{\al}} \Delta_{\CC} \Phi_{\al}
+ g^{ij}\frac{\pa^2f}{\pa w_{\al} \pa w_{\be}}
\frac{\pa \Phi_{\al}}{\pa z^i}\frac{\pa\Phi_{\be}}{\pa z^j} \,.
\end{equation}
Judicious choice of the constants now gives the result, as follows.
First, fix $\ga \in \{ 1, \ldots , n\}$ and choose $C_{\al} = \delta_{\al\ga}$ and $C_{\al \be} = 0$ for all $\al,\be$;
then we deduce that $\Delta_{\CC} \Phi_{\ga} = 0$, giving (\ref{cx-char}a).  Now set $C_{\al} = 0$ for all $\al$ and, for each $\ga = 2, \ldots , n$ in turn, choose $C_{\al \be}$ such that $C_{\al \be} = 0$ for $\al \neq \be$, $C_{\ga \ga} = - C_{11}$, and $C_{\delta\delta}= 0$ for $\delta \neq 1, \ga$.  Then equation (\ref{cx-char}b) follows.  The converse follows from the chain rule \eqref{chain3}.
\end{proof}

\section{Complex-harmonic morphisms and bicomplex manifolds}
\label{sec:cx-ha-bicx}

We now consider the case of maps into a complex surface where we can utilize the
conformal invariance of a complex-harmonic morphism as follows.
By a $2$-dimensional \emph{conformal complex-Riemannian manifold}, we mean a $2$-dimensional complex manifold with an open covering of complex charts on each of which is defined a holomorphic metric, with any two conformally related on the intersection of charts.
By a \emph{complex-orientation} on such a manifold $N$, we mean a reduction of the structure group of $T'N$ to
the positive complex conformal group $\C_+(2,\CC)$.  We claim that a $2$-dimensional complex-oriented conformal complex-Riemannian manifold is the same as a one-dimensional bicomplex manifold.

Indeed, let $N$ be a one-dimensional bicomplex manifold.  In any local bicomplex coordinate $q = z_1 + z_2 \ii_2$, the tensor field
$\dd q \,\dd q^* = \dd z_1^{}\!^2 + \dd z_2^{}\!^2$ defines a holomorphic metric on an open subset of the underlying complex surface,
which we continue to denote by $N$.  Clearly, metrics on overlapping charts are conformally related. In fact, the bicomplex Cauchy-Riemann equations tell us that the structure group of $T'N$ is contained in $\C_+(2,\CC)$ so that $N$ is a $2$-dimensional complex-oriented conformal complex-Riemannian manifold.

Conversely, given such a manifold $N$, locally we can adapt proofs for the semi-Riemannian case such as
\cite[Proposition 14.1.18]{book}, to find complex-oriented complex coordinates $(z_1,z_2)$ which are isothermal in the sense that the conformal class of holomorphic metrics is that containing $\dd z_1{}\!^2 + \dd z_2{}\!^2$.
The transition functions for overlapping charts then have derivatives in $\C_+(2,\CC)$ and so satisfy the bicomplex Cauchy-Riemann equations, hence $N$ has been given the structure of a bicomplex manifold with local bicomplex coordinates $q=z_1 + z_2 \ii_2$\,.

{}From now on we shall consider maps into one-dimensional bicomplex manifolds.  We deduce the following result from Proposition \ref{prop:fund-char}.

\begin{corollary} Let $\Phi : M \to N$ be a holomorphic map from a complex-Riemannian manifold to a one-dimensional bicomplex manifold.  Then $\Phi$ is a com\-plex-har\-mon\-ic morphism if and only if, in any bicomplex chart on $N$,
\begin{equation*} 
 \text{\rm (a)} \quad \Delta_{\CC} \Phi = 0
\quad {\rm and} \quad
\text{\rm (b)} \quad (\grad_{\CC}\Phi)^2 = 0\,.
\end{equation*}
\end{corollary} 

\begin{proof}
In a bicomplex chart, write $\Phi = \Phi_1 + \Phi_2\ii_2$. Clearly, (a) is equivalent to (\ref{cx-char}a).
The equivalence of (b) with (\ref{cx-char}b) follows from the identity
$$
(\grad_{\CC}\Phi)^2
	= (\grad_{\CC}\Phi_1)^2 - (\grad_{\CC}\Phi_2)^2
+ 2\ii_2 \inn{\grad_{\CC}\Phi_1, \grad_{\CC}\Phi_2}_{\CC}\,.
$$
\end{proof}

Note that, if $M$ is an open subset of\/ $\CC^m$,  these equations read
\begin{equation} \label{cx-char2}
\text{\rm (a)} \quad \sum_{k = 1}^m \frac{\pa^2\Phi}{\pa z_k{}\!^2} = 0
\quad {\rm and} \quad \text{\rm (b)} \quad
\sum_{k = 1}^m\left(\frac{\pa\Phi}{\pa z_k}\right)^{\!\! 2} = 0\,.
\end{equation}

Note also that a point is degenerate precisely when $\CN(\grad_{\CC}\Phi) = 0$ but $\grad_{\CC}\Phi \neq \vec{0}$.

Looking at the classification of points in Proposition \ref{prop:types} we see that \emph{any complex-harmonic morphism $\Phi : \CC^m \supset U \to \CC^2 = \BB$ with differential of (complex) rank at most one is degenerate at all points where its differential is non-zero}.

\begin{example} \label{ex:BaWo2}
Embed $\CC$ in $\BB$ as $\CC[i_1]$ by $z \mapsto z+ 0\ii_2$.
Let $U \subset \CC^m$ be open.
Then a smooth map $\Phi:U \to \CC = \CC[\ii_1]$ is a complex-harmonic morphism if and only if it satisfies equations \eqref{cx-char2}
\emph{with $\Phi$ complex-valued}.
Then {\rm (\ref{cx-char2}b)} confirms that $\Phi$ is degenerate away from points where its differential is zero.  This sort of complex-harmonic morphism can be characterized as a map which pulls back holomorphic functions to complex-harmonic ones; for the case $m=4$ see \cite{Ba-Wo-SFR}.
\end{example}

The following proposition gives a way of constructing
com\-plex-har\-mon\-ic morphisms $\Phi$ into a one-dimensional bicomplex manifold $N$ implicitly; it is a bicomplex version of \cite[Theorem 9.2.1]{book}, but care is needed because of the presence of zero divisors.  Since the result is local, by taking a local bicomplex coordinate, we may assume in the proof that $N$ is an open subset of the bicomplex numbers $\BB$; this remark applies to all the results in this section.

\begin{proposition}  \label{prop:impl-B} Let $N$ be a one-dimensional bicomplex manifold and let $A$ be an open subset of\/ $\CC^m \times N$.
Let $\Psi : A \to N$, $(\vec{z},q) \mapsto \Psi(\vec{z},q)$ be a holomorphic function which is bi\-com\-plex-holo\-mor\-phic in its second argument.  Suppose that, for each fixed $q$, the mapping $\vec{z} \mapsto \Psi_q(\vec{z}) := \Psi(\vec{z},q)$ is a complex-harmonic morphism, i.e., satisfies 
\begin{equation} \label{Psi-cond}
\text{\rm (a)}
\quad \sum_{k = 1}^m \frac{\pa^2 \Psi_q}{\pa z_k{}\!^2} = 0
\quad \text{and} \quad
\text{\rm (b)} \quad
\sum_{k = 1}^m \left(\frac{\pa \Psi_q}{\pa z_k}\right)^{\!\! 2} = 0
	\qquad \bigl( (\vec{z},q) \in A \bigr).
\end{equation}

Let $\Phi:U \to N$, $q = \Phi(z)$ be a $C^2$ solution to the equation $\Psi(\vec{z},\Phi(\vec{z})) = \const$ on an open subset\/ $U$ of\/ $M$, and suppose that the mapping $\vec{z} \mapsto \CN(\grad_{\CC}\Psi_q)(\vec{z},\Phi(\vec{z}))$ is not identically zero on $U$.  Then $\Phi$ is a com\-plex-har\-mon\-ic morphism.
\end{proposition}

\begin{proof}  Since
$\vec{z} \mapsto \CN(\grad_{\CC}\Psi_q)(\vec{z},\Phi(\vec{z}))$ is holomorphic but not identically zero, it is non-zero on a dense open subset $\wt{U}$ of $U$.  It suffices to show that $\Phi$ satisfies equations \eqref{cx-char2} on that subset.   {}From the chain rule, at any point $(\vec{z}, \Phi(\vec{z}))$ \ $\bigl(\vec{z} \in \wt{U} \bigr)$ we have
\begin{equation} \label{rule}
\frac{\pa \Psi}{\pa q}\frac{\pa \Phi}{\pa z_i} + \frac{\pa \Psi}{\pa z_i}
 = 0\,.
\end{equation}
Now, at $(\vec{z},\Phi(\vec{z}))$ we have $\CN(\grad_{\CC}\Psi_q) \neq 0$, so that $\CN(\pa \Psi / \pa q) \neq 0$; hence $\pa \Psi/\pa q$ is not a zero divisor. Then, differentiation of $\Psi = 0$ with respect to $\ov{z}_i$ gives $(\pa \Psi / \pa q)(\pa \Phi/\pa \ov{z}_i) = 0$ so that $\pa \Phi/\pa \ov{z}_i = 0$, showing that $\Phi$ is holomorphic.  Again, because $\pa \Psi/\pa q$ is not a zero divisor,  \eqref{rule} gives equation (\ref{cx-char2}b). On differentiating \eqref{rule} once again with respect to
$z_i$, we obtain
\begin{equation} \label{2nd-order-chain}
\frac{\pa \Psi}{\pa q}\frac{\pa^2\Phi}{\pa z_i{}\!^2} +
\frac{\pa^2\Psi}{\pa q^2}\left(\frac{\pa\Phi}{\pa z_i}\right)^{\!\! 2}
+ \frac{\pa^2\Psi}{\pa z_i\pa q}\frac{\pa\Phi}{\pa z_i} 
+ \frac{\pa^2\Psi}{\pa z_i{}\!^2} = 0\,.
\end{equation}
{}From \eqref{rule} we have
$$
\frac{\pa \Psi}{\pa q} \sum_{i = 1}^m \frac{\pa^2\Psi}{\pa z_i\pa q}\frac{\pa \Phi}{\pa z_i} = - \sum_{i = 1}^m \frac{\pa^2\Psi}{\pa z_i\pa q} \frac{\pa \Psi}{\pa z_i} = - \frac{1}{2} \frac{\pa}{\pa q}\sum_{i = 1}^m\left(\frac{\pa \Psi}{\pa z_i}\right)^{\!\! 2} = 0\,.
$$
so that, on summing \eqref{2nd-order-chain} over $i=1, \ldots, m$ and using twice that $\pa \Psi/\pa q$ is not a zero divisor, we obtain equation (\ref{cx-char2}a).  \end{proof} 

\medskip

This leads to a bicomplex version of \cite[Corollary 1.2.4]{book}, with the new feature of degeneracy, as follows.
Write $\vec{\xi}  =\vec{u} + \vec{v}\ii_2$ where
$\vec{u}, \vec{v} \in \CC[\ii_1]^3$. 
  The original case is recovered when $\vec{u}, \vec{v} \in \RR^3$, i.e., $\vec{\xi}$ has values in $\CC^3 = \CC[\ii_2]^3 \subset \BB^3$.  Again, for simplicity, we may assume that $N$ is an open subset of $\BB$.

\begin{corollary} \label{cor:xi}  Let
$N$ be a one-dimensional bicomplex manifold and let
$\vec{\xi} :N \to \BB^3$,
\ $\vec{\xi} = (\xi_1, \xi_2, \xi_3)$, be a
bi\-com\-plex-holo\-mor\-phic map which is \emph{null}, i.e., satisfies
\begin{equation} \label{xi-null}
\vec{\xi}^2 = 0\,,
\end{equation}
and suppose that $\CN(\vec{\xi})$ is not identically zero on $N$.  Then any $C^2$ solution $\Phi: U \to N$, \ $q=\Phi(\vec{z})$, on an open subset of\/ $\CC^3 = \CC[\ii_1]^3$, to the equation
\begin{equation}\label{xieq-B}
\inn{\vec{\xi}(q), \vec{z}}_{\BB} = 1
\end{equation}
 is a com\-plex-har\-mon\-ic morphism of (complex) rank at least one everywhere.  It is degenerate at the points of the fibres $\Phi^{-1}(q)$  \ $(q \in N)$ for which
$\CN\bigl(\vec{\xi}(q)\bigr) = 0$\,.
\end{corollary}

\begin{proof} 
Set
\begin{equation} \label{Psi-xi}
\Psi(\vec{z}, q) = \inn{\vec{\xi}(q)\,,\, \vec{z}}_{\BB}
\qquad (\vec{z} \in \CC^3, \ q \in N)\,.
\end{equation}
Then $\grad \Psi_q = \vec{\xi}(q)$; this is non-zero at any point $q=\Phi(\vec{z})$ by \eqref{xieq-B}.  It follows from Proposition \ref{prop:impl-B} that $\Phi$ is a complex-harmonic morphism; from \eqref{rule} we see that $\dd\Phi \neq 0$ at all points of $U$, so that $\dd\Phi$ has complex rank at least one everywhere.

Let $q \in N$.  On writing $\vec{\xi} = \vec{\xi}(q) = \vec{u} + \vec{v} \ii_2$
where $\vec{u}, \vec{v} \in \CC[\ii_1]^3$, \eqref{xieq-B} is equivalent to the pair of equations
\begin{equation} \label{xieq-B2}
\inn{\vec{u}(q), \vec{z}}_{\CC} = 1, \qquad
\inn{\vec{v}(q), \vec{z}}_{\CC} = 0.
\end{equation}

Note that $\vec{u}$ and $\vec{v}$ span the complex horizontal space
$\Hh^{\cc}_q$ of $\Phi$, and that
$$
\vec{\xi}^2 = \vec{u}^2 - \vec{v}^2
	+ 2\ii_2\inn{\vec{u}, \vec{v}}_{\CC}
\quad \text{and} \quad	
\CN(\vec{\xi})= \vec{u}^2 + \vec{v}^2.
$$
Combining this with \eqref{xi-null} we see that
\begin{equation} \label{xi-2}
\vec{u}^2 = \vec{v}^2 = \frac{1}{2}\CN(\vec{\xi})
	\quad \text{and} \quad \inn{\vec{u},\vec{v}}_{\CC} = 0.
\end{equation}

Suppose that $\CN(\vec{\xi}(q)) \neq 0$.  Then,
$\vec{\xi}(q) \neq \vec{0}$ so that the fibre $\Phi^{-1}(q)$ given by
\eqref{xieq-B} is non-empty.  {}From \eqref{xi-2} we see that
 $\vec{u}$ and $\vec{v}$ are complex-orthogonal with $\vec{u}^2 = \vec{v}^2 \neq \vec{0}$; it follows that they are linearly independent and span a non-degenerate plane.
Hence the fibre \eqref{xieq-B2} is a non-null complex line
which is complex orthogonal to that plane.  By the classification in Proposition \ref{prop:types} (or from \eqref{rule}),
$\Phi$ is submersive at all points on the fibre, with complex horizontal space spanned by $\vec{u}$ and $\vec{v}$.
 
Suppose, instead, that $\CN(\vec{\xi}(q))=0$.
Then from \eqref{xi-null}, $\vec{u}$ and
 $\vec{v}$ span a null subspace of $\CC^3$; since the maximal dimension of such a subspace is one, they  must be linearly dependent.  
Hence, from \eqref{xieq-B2},
the fibre $\Phi^{-1}(q)$ is non-empty if and only if
 $\vec{u} \neq \vec{0}$ but
$\vec{v} = \vec{0}$, in which case it is the
degenerate complex plane
$<\vec{u}(q), \vec{z}>_{\CC} = 1$; from the classification in Proposition \ref{prop:types},
$\Phi$ must be degenerate at each point of this plane.
\end{proof}

We shall now show that any submersive com\-plex-har\-mon\-ic morphism is given locally by Corollary \ref{cor:xi}.

\begin{lemma} \label{lem:lines}  Let $\Phi :U \to N$ be a submersive com\-plex-har\-mon\-ic morphism from an open subset of\/ $\CC^3$ to a one-dimensional bicomplex manifold.  Then the connected components of the fibres of\/ $\Phi$ are open subsets of complex lines in $\CC^3$.
\end{lemma}

\begin{proof}  For convenience, write $\pa_i = \pa/\pa z_i$ \
$\bigl(i=1,2,3, \ (z_1,z_2,z_3) \in U \bigr)$.
Let $p\in U$.  Then, since $\Phi$ is submersive, it is also
 non-degenerate, so we have
$\CN(\grad_{\CC}\Phi )(p) \neq 0$.  Hence we can choose coordinates such that $\pa_1 \Phi(p) = 0$.  Then
\begin{equation} \label{prod-zero}
(\pa_2\Phi + \ii_2\pa_3\Phi) (\pa_2\Phi - \ii_2\pa_3\Phi) = 0
\quad \text{at } p\,.
\end{equation}
Now, since $\CN(\grad_{\CC}\Phi )(p) \neq 0$, one of $(\pa_2\Phi \pm \ii_2\pa_3\Phi)(p)$ must have non-zero complex norm.  Indeed, this follows from the easy calculation at $p$:
\begin{align*}
\CN(\pa_2\Phi + \ii_2\pa_3\Phi) + \CN(\pa_2\Phi - \ii_2\pa_3\Phi)
&= 2\bigl\{\CN(\pa_2\Phi) + \CN(\pa_3\Phi)\bigr\} \\
& = 2\CN(\grad_{\CC}\Phi)\,, \quad \text{since } \pa_1\Phi = 0\,.
\end{align*}
Suppose that $(\pa_2\Phi - \ii_2\pa_3\Phi)(p)$ has non-zero complex norm; the other case is similar.  Then it is not a zero divisor, so from \eqref{prod-zero},  $(\pa_2\Phi + \ii_2\pa_3\Phi)(p) = 0$.  

 On applying the differential operator $\pa_2 - \ii_2\pa_3$ to equation \eqref{prod-zero} and evaluating at $p$, we obtain
$\bigl(\pa_2{}\!^2\Phi(p) + \pa_3{}\!^2\Phi(p)\bigr)
\bigl( \pa_2\Phi(p) - \ii_2 \pa_3\Phi(p)\bigr) = 0$\,, 
so that $\pa_2{}\!^2\Phi(p) + \pa_3{}\!^2\Phi(p) = 0$; then from equation (\ref{cx-char2}a) we obtain $\pa_1{}\!^2\Phi(p) = 0$.

Next, since $p$ is a regular point, we can parametrize the fibre near $p$
by a map $w\to \vec{z}(w) = (z_1(w), z_2(w), z_3(w))$, where each $z_k(w)$ is holomorphic in $w$, and $\vec{z}(0) = p,\ \vec{z}^{\prime}(0) = (1,0,0)$.  Then, by differentiating the equation $\Phi(\vec{z}(w)) = \const$, we obtain
\begin{eqnarray*}
&\sum_{i=1}^3 \pa_i\Phi (\vec{z}(w))\, z_i^{\prime}(w) = 0 \quad \forall\, w\,,\quad \text{and so}
\\ 
&\sum_{i=1}^3\pa_i\pa_j\Phi (\vec{z}(w))\,z_i^{\prime}(w)z_j^{\prime}(w) + \sum_{i=1}^3\pa_i\Phi (\vec{z}(w))\,z_i^{\prime\prime}(w) = 0
\quad \forall\, w\,. 
\end{eqnarray*}
Evaluating the last equation at $w = 0$ gives
$\pa_2\Phi(p)\,z_2^{\prime\prime}(0) + \pa_3\Phi(p)\,z_3^{\prime\prime}(0) = 0$
which can be written as
$\bigl(\pa_2\Phi(p) - \ii_2\pa_3\Phi(p) \bigr)\,(z_2^{\prime\prime}(0) - \ii_2 z_3^{\prime\prime}(0)) = 0$.  
Since $\CN(\pa_2\Phi(p) - \ii_2\pa_3\Phi(p)) \neq 0$, we deduce that $z_2^{\prime\prime}(0) = z_3^{\prime\prime}(0) = 0$.  As the point $p$ was arbitrarily chosen, the lemma follows.
\end{proof}  

To proceed, we make the following assumptions: 
(i) $\Phi:U \to N$ is a submersive complex-harmonic morphism from an open subset of $\CC^3$ onto a one-dimensional bicomplex manifold;
(ii) each fibre component is connected; (iii)  no fibre lies on a complex line through the origin.  Note that,
if (i) holds, after shifting the origin if necessary,
then any point of $U$ has a neighbourhood on which the above conditions are satisfied.

As before, for simplicity in the proof below, we take $N$ to be an open subset of $\BB$.

\begin{proposition} \label{prop:converse} Let $\Phi : U \to N$ be a com\-plex-har\-mon\-ic morphism on an open subset of\/ $\CC^3$ satisfying conditions {\rm (i)--(iii)} above.  Then there is a unique bi\-com\-plex-holo\-mor\-phic map $\vec{\xi} : N \to \BB^3$ with $\vec{\xi}^2 = 0$ and $\CN(\vec{\xi} ) \neq 0$ such that the fibre of\/ $\Phi$ at $q \in N$ is given by \eqref{xieq-B}.
\end{proposition}

\begin{proof}  Let $\ell_0 = \ell_0(q)$ be the complex line through the origin parallel to $\Phi^{-1}(q)$ and set $\Pi = \ell_0^{\bot_{\CC}}$ $ := \{ \vec{w}\in \CC^3 : \inn{\vec{z}, \vec{w}}_{\CC} = 0 \ {\rm for \ all} \ \vec{z} \in \ell_0\}$.  Since $\Phi$ is submersive, $\ell_0$ is not null so that $\Pi \cap \ell_0$ is a single point, $\vec{c}$, say, with $\vec{c}^2 \neq 0$.

Recalling that $\grad\Phi_1$ and $\grad\Phi_2$ are complex-orthogonal with the same non-zero complex norm, set 
$$\vec{\ga} = \grad \Phi_1 \times \grad\Phi_2 /(\grad\Phi_1)^2 =
	 \ii_2\,\grad \Phi \times \grad\Phi^* /\CN(\Phi);
$$
then $\vec{\ga}$ is one of the two vectors of complex norm $1$ parallel to $\ell_0$.  
Set $J\vec{c} = \vec{\ga} \times_{\CC} \vec{c}$ where $\times_{\CC}$ denotes the vector product in $\RR^3$ extended to $\CC^3$ by complex bilinearity.  Next set
 $\vec{\xi} = \vec{\xi}(q)
= (\vec{c} + \ii_2 J\vec{c})\big/ \vec{c}^2$,  so that $\vec{\xi}^2 = 0$. Then the fibre is given by $\inn{\vec{\xi}(q), \vec{z}}_{\BB} = 1$;
since $\vec{c}^2 \neq 0$, this is well-defined, and $\CN(\vec{\xi}) \neq 0$.

It remains to prove that $\vec{\xi}:N \to \BB^3$ is
bi\-com\-plex-holo\-mor\-phic.  To do this, we show that
$\pa \vec{\xi}/\pa q^*=0$ in a way analogous to
\cite[Lemma 1.3.3]{book}.

Let $q_0\in N$ and let $\vec{z}^0 \in \Phi^{-1}(q_0)$.
In the following calculations, all quantities are evaluated at $\vec{z}^0$ or $q_0$.
As in Lemma \ref{lem:lines}, we may suppose that our coordinates are chosen such that, at $z^0$,
\begin{equation}
\pa_1\Phi = 0\,, \ 
\pa_2\Phi + \ii_2 \pa_3\Phi = 0 \ \text{and} \
\CN(\pa_2\Phi - \ii_2 \pa_3\Phi) \neq 0\,.
\label{phieq0b}
\end{equation}
Further, without loss of generality, we may choose the coordinates so that
$\vec{z}^0$ is the point $(0,0,1)$.  Then the fibre $\Phi^{-1}(q_0)$ through
$\vec{z}^0$ is a connected open subset of the complex line parametrized by
$w \mapsto \vec{z}(w) = (w,0,1)$.  
On applying the operator $\pa_2 + \ii_2 \pa_3$ to equation
(\ref{xieq-B}) we obtain
\begin{equation}
\Biginn{\,\frac{\pa\vec{\xi}}{\pa q}\, \bigl(\pa_2\Phi
	+ \ii_2 \pa_3\Phi \bigr) +
\frac{\pa\vec{\xi}}{\pa q^*}\,
	\bigl(\pa_2 \Phi^* + \ii_2
		\pa_3 \Phi^*\bigr) \,, \ \vec{z}(w)}_{\!\BB}  \,
	+ \, \xi_2 + \xi_3 \ii_2= 0 \,.
\label{xidiff} \end{equation}
Now $\CN(\pa_2\Phi^* + \ii_2 \pa_3\Phi^*)
= \CN(\pa_2\Phi - \ii_2 \pa_3\Phi) \neq 0$ at $\vec{z}^0$. 
By continuity and connectedness of the fibres, (\ref{phieq0b}) holds at all points of the fibre.
Also,  on the fibre we have
\begin{equation}
\xi_1   =  0  \quad \text{and} \quad
	(\xi_2 + \xi_3\ii_2)(\xi_2 - \xi_3\ii_2)   =  0   \,.
\label{xi3} \end{equation}
Now at $\vec{z}^0$, if we write $\grad\Phi_1 = (0,a,b)$, then  
$\grad\Phi_2 = \pm(0,-b,a)$.  With the minus sign, this gives
$\pa_2\Phi + \ii_2\pa_3\Phi = 0$ 
in contradiction to \eqref{phieq0b}, hence $\grad\Phi_2 = +(0,-b,a)$
and we have 
$ \vec{\ga} = (0,a,b) \times (0,-b,a) \big/(a^2+b^2) = (1,0,0)$.
Since $\vec{c} = (0,0,1)$, this gives
$J\vec{c} = (0,-1,0)$ and
$\vec{\xi}(q_0) = (0,-\ii_2,1)$,
so that $\xi_2 - \xi_3\ii_2$ is not a divisor of zero, and from \eqref{xi3}  we see that 
$\xi_2 + \xi_3\ii_2 = 0$ on the fibre.
Then (\ref{xidiff}) becomes
$$
\frac{\pa\xi_1}{\pa q^*}\, w + \frac{\pa\xi_3}{\pa q^*} = 0
\,.
$$
Since this is valid for all $w$ in a neighbourhood of $0$, we conclude that
$$
\frac{\pa\xi_1}{\pa q^*} = \frac{\pa\xi_3}{\pa q^*} = 0 \,.
$$
On the other hand, on differentiating $\vec{\xi}^2 = 0$ and evaluating at $q_0$ we obtain
$\xi_2 \,(\pa \xi_2 / \pa q^*) = 0$. Now  $\xi_2 = -\ii_2$ is not a zero divisor; so we conclude that
$$
\frac{\pa\xi_1}{\pa q^*} = \frac{\pa\xi_2}{\pa q^*} =
\frac{\pa\xi_3}{\pa q^*} = 0
$$
at $\vec{z}^0$.  Since $\vec{z}^0$ is an arbitrary point of $N$, this shows that $\vec{\xi}$ is bicomplex-holomorphic.
\end{proof}

\begin{remark} \label{rem:orientation}
{\rm (i)} We see that $\vec{\ga}$ gives the direction of the fibres, oriented as explained below, and $\vec{c}$ gives their displacement from the origin; we call $\vec{\ga}$ and $\vec{c}$ the \emph{Gauss map} and \emph{fibre position map} of\/ $\Phi$, respectively.  

{\rm (ii)}The process of picking one of the two possible values of $\vec{\ga}$ may be explained as follows.
Let $\Pi$ be a non-degenerate complex $2$-plane in $\CC^3$ and let $\vec{u}, \vec{v}$ be a complex-orthogonal basis with
$\vec{u}^2 = \vec{v}^2 \ (\neq 0)$.  A (complex-) orientation of\/ $\Pi$ is an equivalence class of such bases under the equivalence relation that they are related by a member of $\C_+(2,\CC)$ (see \eqref{C+}).  In particular, two complex-orthonormal bases have the same orientation if and only if they are related by a member of $SO(2,\CC)$.  To any complex-oriented plane, there is a canonical complex-normal of complex norm one, given by $\vec{u} \times \vec{v}/\vec{u}^2$ for any oriented complex-orthogonal basis with $\vec{u}^2 = \vec{v}^2$; call it the \emph{oriented complex-normal}.

In the above proof, we are lifting the canonical complex-orientation of the codomain to a complex-orientation of the complex-horizontal space,  and then $\vec{\ga}$ is its oriented complex-normal.
\end{remark}

We can find all triples $\vec{\xi} = (\xi_1, \xi_2, \xi_3)$ of bi\-com\-plex-holo\-mor\-phic functions satisfying $\vec{\xi}^2 = 0$, i.e., $\sum_k\xi_k{}\!^2 = 0$, as in the complex case.  Indeed, provided that
$\xi_2 - \xi_3\ii_2$ is not a zero divisor, there are bi\-com\-plex-holo\-mor\-phic functions $G$ and $H$ with $\CN(H) \neq 0$ such that
\begin{equation} \label{bicomplex-Weier}
(\xi_1, \xi_2, \xi_3)
= \frac{1}{2H}\bigl(-2G,\, 1-G^2,\, (1+G^2)\ii_2\bigr).
\end{equation}
To see this, as for the Riemannian Weierstrass representation, it suffices to take $G = -\xi_1 \big/ (\xi_2 - \xi_3 \ii_2)$ and
$H = 1 \big/ (\xi_2 - \xi_3 \ii_2)$\,.

The equation $\inn{\vec{\xi}(q), \vec{z}}_{\BB} = 1$ then reads
\begin{equation} \label{W-cx}
 -2G\,z_1 + (1-G^2)z_2 + (1+G^2)z_3\ii_2 = 2H\,;
 \end{equation}
 note that, in contrast to \eqref{bicomplex-Weier}, this makes sense even when $H=0$.

 \section{Interpretation and Compactification}
 \label{sec:interp}

Given bicomplex-holomorphic functions $q \mapsto G(q)$ and
$q \mapsto H(q)$ defined on an open subset $N$ of $\BB$, or more generally on a one-dimensional bicomplex manifold, we can form the equation
\eqref{W-cx}. By Corollary \ref{cor:xi}, $C^2$ solutions $q = \Phi(\vec{z})$ to this equation are complex-harmonic morphisms from open subsets of $\CC^3$ to $N$, and by Proposition \ref{prop:converse}, all such harmonic morphisms which are submersive are given this way, locally.  In general,
the equation \eqref{W-cx} defines a \emph{congruence of lines and planes}; indeed, for each $q \in N$, if $\CN(G) \neq -1$,  \eqref{W-cx} defines a complex line, whereas if $\CN(G) = -1$, \eqref{W-cx} either has no solutions or defines a plane, see Proposition \ref{prop:parallel} below. We shall call these lines and planes the \emph{fibres} of the congruence
\eqref{W-cx} as they form the fibres of any smooth harmonic morphism $q = \Phi(\vec{z})$ which satisfies that equation. However, starting with arbitrary data $G$ and $H$, the fibres of the congruence
\eqref{W-cx} may intersect or have \emph{envelope points} where they become infinitesimally close. We shall consider the behaviour of this congruence when the fibres are degenerate or have direction not represented by a finite value of $G$.
We consider first non-degenerate fibres.

Recall the standard chart of $S^2_{\CC}$ given by complexified stereographic projection \eqref{stereo}.
Then, as in \cite{Ba-Wo-Bernstein}, it is easy to see that $\vec{\ga} = \si_{\CC}^{-1}G$ is the Gauss map giving the oriented direction of the fibre and
$\vec{c} = (\dd\si_{\CC}^{-1})_G(H)$ is the fibre position map, as defined in Remark \ref{rem:orientation}.

Let $\CP^2$ denote complex projective $2$-space and let $\Zz = \{[z_1,z_2,z_3] \in \CP^2: z_1^{}\!^2 + z_2^{}\!^2 +z_3^{}\!^2 = 0 \}$; thus points of $\Zz$ represent null one-dimensional complex subspaces of $\CC^3$.  We have a $2\!:\!1$ mapping $S^2_{\CC} \to \CP^2 \setminus \Zz$ given by $\vec{z} \mapsto [\vec{z}]$; the image of $\vec{\ga}$ under this mapping is the complex line through the origin which is parallel to the fibre, with its complex-orientation forgotten. 
 
An alternative interpretation is as follows.  Let
\begin{equation} \label{CQ1B*}
\CQ^1_{\BB *} = \{ \vec{\xi} = (\xi_1,\xi_2,\xi_3) \in \BB^3 : \vec{\xi}^2 = 0, \ \CN(\vec{\xi}) \neq 0\}.
\end{equation}
 For $\vec{\xi} \in \CQ^1_{\BB *}$, write 
$\vec{\xi} = \vec{u} + \vec{v}\ii_2$ with
$\vec{u}, \vec{v} \in \CC^3$.  Then
$\vec{u}^2 = \vec{v}^2 = \frac{1}{2}\CN(\vec{\xi}) \neq 0$ and $\inn{\vec{u},\vec{v}}_{\CC} = 0$.  
Projectivizing $\CQ^1_{\BB *}$ gives the open dense subset
$\Q^1_{\BB *}$ of the bicomplex quadric given by \eqref{Q1B*}.  
Let $G_2(\CC^3)$ be the Grassmannian of $2$-dimensional complex subspaces in $\CC^3$ and let $\Dd$ denote the set of points in $G_2(\CC^3)$ which represent degenerate $2$-dimensional subspaces.
Note that the condition $\CN(\vec{\xi}) \neq 0$ is equivalent to linear independence of the vectors $\vec{u}$ and $\vec{v}$ so that they span a complex $2$-dimensional subspace; hence we have a double covering $\Q^1_{\BB*} \to G_2(\CC^3) \setminus \Dd$ given by
$[\vec{\xi}] = [\vec{u} \pm \vec{v}\ii_2] \mapsto \spn\{\vec{u},\vec{v}\}$, thus we can think of $\Q^1_{\BB*}$ as the space of \emph{complex-oriented} non-degenerate
$2$-dimensional subspaces of $\CC^3$. 

Now we have a map $\Q^1_{\BB *}  \to S^2_{\CC}$ given by 
\begin{equation} \label{Q1B*-S2C}
[\vec{\xi}]= [\vec{u} + \vec{v}\ii_2]
	\mapsto \vec{u}\times\vec{v}/\vec{u}^2 = 
		\vec{u}\times\vec{v}/\vec{v}^2
			= (\vec{\xi} \times \vec{\xi}^*)\ii_2  \big/\CN(\vec{\xi})\,.
\end{equation}
This is well-defined and is an equivalence of bicomplex manifolds; see Remark \ref{rem:maps} for its inverse.
Furthermore, it covers the biholomorphic map
$G_2(\CC^3)\setminus \Dd \to \CP^2 \setminus \Zz$ given by sending a subspace $\spn\{\vec{u},\vec{v}\}$ to its orthogonal complement
$[\vec{u} \times \vec{v}]$.

We have thus established the bottom left-hand square of the commutative diagram below in which all spaces are
two-dimensional conformal complex-Riemannian manifolds and all maps between them are holomorphic.
Further, all three spaces in the middle row are 
one-dimensional \emph{bi}complex manifolds and the top three vertical arrows are the standard charts of Examples \ref{ex:S2C}--\ref{ex:complex-quadric}.  
The maps in the first commutative diagram are as shown in the second diagram where, for brevity, we write $C = \CN(G)$.

\begin{diagram}[tight,height=0.6cm,width=2.2cm]
 \BB \setminus \Hh^1 &  \rTo^{\rm Id} & \BB \setminus \Hh^1 & \rInto^{\rm inclusion} & \BB\\
\dTo & & \dTo & & \dTo \\
\Q^1_{\BB *} & \rTo^{\cong} & S^2_{\CC}  & \rInto^{\space\iota_{S^2_{\CC}}} & \Q^2_{\CC} \\
\dTo^{\text{2:1}} & & \dTo_{\text{2:1}} & &
\dTo_{\text{branched 2:1}} \\ 
G_2(\CC^3) \setminus \Dd 	& \rTo^{\cong} & \CP^2 \setminus \Zz & \rInto & \CP^2
\end{diagram}

\bigskip

\small
\begin{diagram}[height=0.6cm,width=2.2cm]
 G =G_1+G_2\ii_2 &  \rMapsto & G =G_1+G_2\ii_2&
	\rMapsto& G =G_1+G_2\ii_2\\
\dMapsto & & \dMapsto & & \dMapsto \\
 [-2G, 1-G^2, (1+G^2)\ii_2]   & \rMapsto &
	(1-C, 2G_1, 2G_2)/(1+C) & \rMapsto &
		[1+C,  1-C, 2G_1, 2G_2]  \\
= & & = & & = \\	
\vec{\xi} = \vec{u} + \vec{v}\ii_2  && \vec{u} \times \vec{v}/\vec{u}^2 
		& & [\vec{u}^2, \vec{u} \times \vec{v} ]  \\
\dMapsto & & \dMapsto & & \dMapsto \\ 
 \spn(\vec{u}, \vec{v}) & \rMapsto & [\vec{u}\times\vec{v}]&
 	\rMapsto & [\vec{u} \times \vec{v}]
\end{diagram}

\normalsize \bigskip

The Gauss map $\ga$ is a map from $N$ to $\Q^1_{\BB *}$ or, equivalently, $S^2_{\CC}$\,.
The fibre position map $\vec{c}$ is a map from $N$ to the tautological bundle
$\CQ^1_{\BB *} \to \Q^1_{\BB *}$, see \eqref{CQ1B*},
or to the holomorphic tangent bundle of $S^2_{\CC}$, which covers $\ga$. 

In order to include degenerate fibres and directions corresponding to values of $G$ `at infinity', we compactify this picture as follows.  There is a natural  bicomplex-holomorphic inclusion map $\iota_{S^2_{\CC}}: S^2_{\CC} \hookrightarrow \Q^2_{\CC}$ defined by $(\zeta_1,\zeta_1,\zeta_3) \mapsto [1,\zeta_1,\zeta_2,\zeta_3]$ (see Example \ref{ex:complex-quadric}).   In the standard charts of Examples \ref{ex:S2C} and \ref{ex:complex-quadric}, this is given by
$G \mapsto [1+G_1{}\!^2+G_2{}\!^2, 1-G_1{}\!^2-G_2{}\!^2,2G_1,2G_2]$.

The double cover $S^2_{\CC} \to \CP^2 \setminus \Zz$ extends to a map $\Q^2_{\CC} \to \CP^2$ given by forgetting the first component. This is surjective, and is $2\!:\!1$ away from $\Zz$ where it is branched.  

Degenerate fibres appear if we allow  $\CN(\vec{\xi}) = 0$, i.e., $[\vec{\xi}] \in \Q^1_{\BB}\setminus \Q^1_{\BB *}$; in the standard chart for $\Q^1_{\BB}$, this corresponds to $\CN(G) = -1$. Then $\vec{u}$ and $\vec{v}$ become collinear null complex vectors, and the horizontal space, $\spn\{\vec{u}, \vec{v}\}$, collapses to a null complex \emph{line}.  Its complex-orthogonal complement is a degenerate complex \emph{plane} through the origin; 
the fibre is either empty or a degenerate complex plane parallel to this (see Proposition \ref{prop:parallel} below).  We get no point in $S^2_{\CC}$ but we do get points in $\Q^2_{\CC}$, and thus in $\CP^2$, as explained by the following two lemmas.  Recall the fattened origin
$\Nn = \{\vec{\xi}\in \BB^3 : \CN(\xi_i)=0\ \forall i\}$.  The following explains our earlier use of $\Nn$.

\begin{lemma} \label{lem:complex-rep}  Let $\vec{\xi} \in \BB^3 \setminus \Nn$ have $\vec{\xi}^2 = 0$.  Then

{\rm (i)} $\CN(\vec{\xi}) = 0$
if and only if

{\rm (ii)} there exists $\vec{\xi}_{\CC} \in \CC[\ii_1]^3 \setminus \{ \vec{0}\}$ with $\vec{\xi}_{\CC}{}^2 = 0$ such that  $\vec{\xi} = \la \vec{\xi}_{\CC}$ for some $\la \in \BB$
with $\CN(\la )\neq 0$.

Further, if condition {\rm (ii)} holds, the projective class $[\vec{\xi}_{\CC}]\in \CP^2$ of\/  $\vec{\xi}_{\CC}$ is unique.
\end{lemma}

\begin{proof}  If (ii) holds, then
$\CN(\vec{\xi}_i) = \CN(\lambda)\,(\vec{\xi}_{\CC})_i{}^2$.  Since
$\vec{\xi}_{\CC} \neq \vec{0}$, we have
$\CN(\vec{\xi}_i) \neq 0$ for some $i$, i.e., $\vec{\xi} \not\in \Nn$.  Also
$\CN(\vec{\xi}) = \CN(\lambda)\,\vec{\xi}_{\CC}{}^2 = 0$, so (i) holds.

Conversely, if (i) holds, write
$\vec{\xi} = \vec{u} + \vec{v} \ii_2$ with
$\vec{u}, \vec{v} \in \CC[\ii_1]$\,.  Then since both $\vec{\xi}^2 = 0$ and $\CN(\vec{\xi}) = 0$, we have $\vec{u}^2 = \vec{v}^2 = \inn{\vec{u}, \vec{v}}_{\CC} = 0$.  Then \emph{either} $\vec{u} \neq 0$ and $\vec{v} = \mu \vec{u}$ for some $\mu \in \CC$, \emph{or} $\vec{v} \neq 0$ and $\vec{u} = \nu \vec{v}$ for some $\nu \in \CC$.
 
Without loss of generality, we may assume that we have the first case.  Then $\vec{\xi} = \la\vec{u}$ where $\la = 1 + \mu \ii_2$ and we set $\vec{\xi}_{\CC} = \vec{u}$.  
Since $\vec{\xi} \not\in \Nn$, we have $\CN(\la ) \neq 0$, so (ii) holds. 

For uniqueness of $[\vec{\xi}_{\CC}] \in \CP^2$, given two representations
$\vec{\xi} = \la \vec{\xi}_{\CC} = \la^{\prime} \vec{\xi}_{\CC}{}^{\prime}$, then 
$\vec{\xi}_{\CC}{}^{\prime}
= (\la^{\prime})^{-1}\la\, \vec{\xi}_{\CC}$
and necessarily
$(\la^{\prime})^{-1} \la \in \CC \setminus \{0\}$.
\end{proof}

We shall call $\vec{\xi}_{\CC}$ a \emph{complex representative of\/ $\vec{\xi}$}, and $[\vec{\xi}_{\CC}]$ its \emph{complex projective representative}. 

\begin{proposition}
The bicomplex-holomorphic diffeomorphism
$\Q^1_{\BB *} \to S^2_{\CC}$ extends to a bicomplex-holomorphic diffeomorphism $\phi : \Q^1_{\BB} \to \Q^2_{\CC}$ given by
\begin{equation} \label{Q1B-to-Q2C}
\phi ([\vec{\xi}]) = \left\{ \begin{array}{ll}
\bigl[\CN(\vec{\xi}),\, (\vec{\xi} \times \vec{\xi}^*)\ii_2 \bigr]
	& \qquad 	\bigl(\CN(\vec{\xi}) \neq 0 \bigr),\\
\text{$ \bigl[0, \vec{\xi}_{\CC}\bigr]$}
 & \qquad \bigl(\CN(\vec{\xi}) = 0 \bigr),
\end{array} \right.  
\end{equation}  
where $\vec{\xi}_{\CC}$ is a complex representative of $\vec{\xi}$.
\end{proposition}
 \begin{proof}
 Note that the first formula applies to points of $\Q^1_{\BB *}$ and the second to points of $\Q^1_{\BB} \setminus \Q^1_{\BB *}$.
 
 First we show that the map $\phi$ is well-defined.  If $[\vec{\eta}] = [\vec{\xi} ]$ then $\vec{\eta} = \la \vec{\xi}$ with $\CN(\la ) \neq 0)$
so that $\CN(\vec{\eta}) = \CN(\la)\CN(\vec{\xi})$ and
$\vec{\eta} \times \vec{\eta}^* = \la \la^*\vec{\xi} \times \vec{\xi}^* = \CN(\la ) \vec{\xi}\times \vec{\xi}^*$. 
Hence, if $\CN(\vec{\xi}) \neq 0$, then
$[\CN(\vec{\eta}),\, (\vec{\eta} \times \vec{\eta}^*)\ii_2]
= [\CN(\vec{\xi}),\, (\vec{\xi} \times \vec{\xi}^*)\ii_2]$. 
On the other hand, if $\CN(\vec{\xi}) = 0$, then
$\phi([\vec{\xi}] = [0,\vec{\xi}_{\CC}]$, which is well-defined by uniqueness of
$[\vec{\xi}_{\CC}]$.

Now, in the standard chart for $\Q^1_{\BB_*}$, the map $\phi$ is given by
$$
\phi([\vec{\xi}]) = [1+\CN(G), 1-\CN(G), 2 G_1, 2 G_2] \,,
$$
with similar expressions in the other charts for $\Q^1_{\BB_*}$.
This shows that $\phi$ is smooth, in fact complex analytic; to see that it is bicomplex-holomorphic, note that, in the standard chart for $\Q^2_{\CC}$, it is just the identity map $G \mapsto G$, and similarly in the other charts.

In order to prove that $\phi$ is a diffeomorphism, we need to find a (two-sided) smooth inverse $\psi$.  Using the charts $G$, $\check{G}$, $L$
and $K$ for $\Q^2_{\CC}$ (Example \ref{ex:complex-quadric}), we obtain 
\begin{multline*}
\psi([\vec{\zeta}]) =
   \bigl[- 2(\zeta_0+\zeta_1)(\zeta_2+ \zeta_3\ii_2),
	(\zeta_0+ \zeta_1)^2- (\zeta_2 + \zeta_3\ii_2)^2, \\[-1.2ex]
		\bigl( (\zeta_0+\zeta_1)^2 + ( \zeta_2 + \zeta_3 \ii_2)^2 \bigr) \ii_2
	\bigr] \qquad ([\vec{\zeta}] \in V_G),
	\end{multline*}	\\[-9ex]
\begin{multline*} 
\psi([\vec{\zeta}]) =
 \bigl[- 2(\zeta_0-\zeta_1)(\zeta_2- \zeta_3\ii_2),
	 -(\zeta_0- \zeta_1)^2+ (\zeta_2 - \zeta_3\ii_2)^2,\\[-1.3ex]
	\bigl( (\zeta_0-\zeta_1)^2 + ( \zeta_2 - \zeta_3 \ii_2)^2 \bigr) \ii_2
	\bigr] \qquad ([\vec{\zeta}] \in V_{\check{G}}),
\end{multline*}  \\[-9ex]
\begin{multline*}
\psi([\vec{\zeta}]) = \bigl[\bigl( (\zeta_0+\zeta_2)^2 + ( \zeta_3 + \zeta_1\ii_2)^2\bigr) \ii_2, - 2(\zeta_0+\zeta_2)(\zeta_3 + \zeta_1\ii_2),     \\[-1.2ex]
		(\zeta_0 + \zeta_2)^2-(\zeta_3+\zeta_1\ii_2)^2
	\bigr] \qquad ([\vec{\zeta}] \in V_L),
	\end{multline*}	\\[-9ex]
\begin{multline*} 
\psi([\vec{\zeta}]) =
 \bigl[  (\zeta_0 + \zeta_3)^2 - ( \zeta_1 + \zeta_2\ii_2)^2,
	 \bigl( (\zeta_0 + \zeta_3)^2 + ( \zeta_1 + \zeta_2 \ii_2)^2 \bigr) \ii_2, \\[-1.2ex]
  - 2( \zeta_0 + \zeta_3)( \zeta_1 + \zeta_2 \ii_2) \bigr]
   \qquad ([\vec{\zeta}] \in V_K).
   \end{multline*}
   
 That on the intersections of charts, the above expressions for $\psi$ coincide is readily checked using the defining relation $-\zeta_0{}^2 + \zeta_1{}^2 + \zeta_2{}^2 + \zeta_3{}^2 = 0$ of $\Q^2_{\CC}$.  The map $\psi$ is clearly complex analytic and it can be checked that it really is a two-sided inverse for $\phi$, so is bicomplex-holomorphic; we omit the calculations.
\end{proof}

\begin{remark} \label{rem:maps}
{\rm (i)} The restriction of $\psi$ to $S^2_{\CC}$
gives an inverse to the map \eqref{Q1B*-S2C}.

{\rm (ii)} The map \eqref{Q1B-to-Q2C} restricts to a bijection from the set of
null directions
$\Q^1_{\BB} \setminus \Q^1_{\BB*}$ to the points at infinity
$\{[0,\zeta_1,\zeta_2,\zeta_3] \in \Q^2_{\CC}\} \cong \CP^1$; thus $S^2_{\CC}$ has been conformally compactified by adding a $\CP^1$, giving the compactification $\Q^1_{\BB}$, or equivalently $\Q^2_{\CC}$\,.

{\rm (iii)} The branched double cover $\Q^2_{\CC} \to \CP^2$ is one-to-one on the points at infinity of $\Q^2_{\CC}$, and maps  $[0,\zeta_1,\zeta_2,\zeta_3]$ to the point $[\zeta_1,\zeta_2,\zeta_3]$ of $\Zz$.

{\rm (iv)} In the standard chart for $\Q^1_{\BB}$ and $\Q^2_{\CC}$, the direction
$[\vec{\xi}]$ is null when $\CN(G) = -1$; then $\phi([\vec{\xi}])$ is the point at infinity $[0,1,G_1,G_2] \in \Q^2_{\CC}$.  The double cover
$\Q^2_{\CC} \to \CP^2$ maps this to $[1,G_1,G_2] \in \Zz$.
\end{remark}

Now note that the double cover $\Q^1_{\BB *} \to G_2(\CC^3) \setminus \Dd$, \ 
$[\vec{\xi}] = [\vec{u} + \vec{v}\ii_2] \mapsto \spn\{\vec{u}, \vec{v}\}$,
extends to a double cover $\Q^1_{\BB} \to G_2(\CC^3)$ given on 
$\Q^1_{\BB} \setminus \Q^1_{\BB *}$ by $[\vec{\xi}] \mapsto [\vec{\xi}_{\CC}]^{\perp_{\CC}}$, where
$[\vec{\xi}_{\CC}]$ is the complex projective representative of $[\vec{\xi}]$ as defined in Lemma \ref{lem:complex-rep}.  That this is holomorphic is easily checked.

We thus obtain the following commutative diagram which extends the previous commutative diagram above to include degenerate directions; all maps are bicomplex-holomorphic.  

\begin{diagram}[tight,height=0.6cm,width=2.2cm]
 \BB  &\rTo^{\rm Id} & \BB\\
\dTo & &  \dTo \\
\Q^1_{\BB} & \rTo^{\cong} & \Q^2_{\CC} \\
\dTo^{\text{branched double cover}} & &
\dTo_{\text{branched double cover}} \\ 
G_2(\CC^3) & \rTo^{\perp_{\CC}} & \CP^2
\end{diagram}

\bigskip

Finally, the behaviour of $H$ at a degenerate fibre is described by the following result.

\begin{proposition} \label{prop:parallel}
{\rm (i)} Suppose that $\CN(G) \neq -1$.  Then the equation  \eqref{W-cx} represents a non-null line.

{\rm (ii)} Suppose that $\CN(G) = -1$.  Then the equation
\eqref{W-cx} has solutions if and only if\/ $H$ is a complex multiple of\/ $G$, in which case it represents a degenerate plane.
\end{proposition}

\begin{proof} 
Writing $G = G_1 + G_2\ii_2,\, H = H_1 + H_2 \ii_2$, the equation \eqref{W-cx} is equivalent to the pair of complex equations
$$
\left\{ \begin{array}{rcl}
 - 2G_1 z_1 + (1-G_1{}^2 + G_2{}^2)z_2 - 2G_1G_2 z_3 & = & 2H_1\,,\\
 - 2G_2 z_1 - 2G_1G_2z_2 + (1 + G_1{}^2 - G_2{}^2)z_3 & = & 2H_2\,.
\end{array}\right.
$$
This defines a line
unless the left-hand side coefficients of the two equations are proportional,
 which happens precisely when $\CN(G) = -1$.  In that case, the pair becomes
$$
G_1(z_1 + G_1z_2 + G_2 z_3) = - H_1\,, \qquad
	G_2(z_1 + G_1 z_2 + G_2z_3) = -H_2\,;
$$
this has a solution if and only if $H$ is a complex multiple of $G$, in which case it reduces to one equation  and so defines a plane.  This plane is easily seen to be degenerate, indeed the vector $[1,G_1,G_2]$  is both complex-normal and parallel to it.
\end{proof}

\begin{remark} On using the formula
$\vec{c} = (\dd\si_{\CC}^{-1})_G(H)$, we can easily show that, as we approach a degenerate fibre, the fibre position map $\vec{c}$ becomes collinear with $\vec{\ga}$
and has complex norm which grows as $1/\CN(\xi)$.
\end{remark}

\begin{example} \label{ex:proj-cx}
\rm (Complex orthogonal projection)

Put $G=0$, $H=(1/2)q$. Then equation \eqref{W-cx} becomes
$$
z_2 + z_3\ii_2  = q
$$
which has solution $q = \phi(\vec{z}) = z_2 + z_3\ii_2$.  This is simply an orthogonal projection $\CC^3 \to \CC^2$.
\end{example}

\begin{example} \label{ex:radial-cx}
\rm (Complex radial projection)
Put $G = q$, $H=0$, then \eqref{W-cx} becomes
the quadratic equation
\begin{equation} \label{W-radial-cx}
(z_2-z_3\ii_2)q^2 + 2 z_1 q - (z_2+z_3\ii_2) = 0\,.
\end{equation}

Let $U$ be an open set in $\CC^3 \setminus \{z_2=z_3 = 0\} \setminus \{z_1{}\!^2 + z_2{}\!^2 + z_3{}\!^2 = 0 \}$ on which there is a smooth branch of $\sqrt{z_1{}\!^2 + z_2{}\!^2 + z_3{}\!^2\,}$, then \eqref{W-radial-cx} has four solutions $q(\vec{z})$ with $z \in U$:
\begin{equation} \label{radial-cx}
q = \bigl(-z_1 + \ve\sqrt{z_1{}\!^2+z_2{}\!^2+z_3{}\!^2}\,\bigr)
\big/ (z_2-z_3\ii_2) \qquad (\ve = \pm 1, \pm\jj)\,.
\end{equation} 

When $\ve = \pm 1$, $q = \si_{\CC}\bigl(\pm z / 
	\sqrt{z_1{}\!^2 + z_2{}\!^2 + z_3{}\!^2}\,\bigr)$,
i.e., it is the complexification of $\pm$ radial projection
$\RR^3 \setminus \{\vec{0}\} \to S^2$ composed with stereographic projection (see \cite[Example 1.5.2]{book}).

When $\ve = \pm j$, we have $qq^* = -1$, so that \eqref{radial-cx} defines an everywhere-degenerate harmonic morphism with fibres the complex $2$-planes tangent to the light cone
$z_1{}\!^2 + z_2{}\!^2 + z_3{}\!^2 = 0$.

For comparison with the semi-Riemannian cases below, note that
$G = q\ii_1$, $H=0$ gives the same map up to the isometry $q \mapsto q\ii_1$.
\end{example} 

\begin{example} \label{ex:disc-cx} \rm (Complex disc example)
Put $G(q) = q$ and $H(q) = t\,q\,\ii_2 $ where $t \in \CC[\ii_1]$ is a complex number.  Then \eqref{W-cx} becomes the quadratic equation
\eqref{W-radial-cx} with $z_1$ replaced by $z_1 + t\ii_2$.
This again has four solutions
$\vec{z} \mapsto q(\vec{z})$ on suitable domains.  

Again, note that
$G = q\ii_1$, $H= tq\ii_1\ii_2 = t\jj z$ gives the same map up to the isometry $q \mapsto q\ii_1$.
\end{example}

\begin{remark} \rm
There are many complex-harmonic morphisms from open subsets of $\CC^3$ to $\CC^2 = \BB$ which are not obtained by extending a real harmonic morphism.  Indeed, as in Remark \ref{rem:Ringleb}, write $q = z a + w b$ and  take
$G(q) = g_1(z) a + g_2(w) b$ and $H = h_1(z)a+h_2(w)b$.  Then if $\Phi$ is the extension of a harmonic morphism on a domain of $\RR^3$, we must have $g_1 = g_2$ and $h_1=h_2$.
\end{remark}

\section{Real harmonic morphisms}\label{sec:real}

Harmonic morphisms from open subsets of $\RR^3$ to $\RR^2$
were discussed in \cite{Ba-Wo-Bernstein} and
\cite[Chapter 1]{book}; they are recovered
 from our theory by setting $z_i$ real, taking $\Phi$ with 
values in $\CC$, and embedding $\CC$ in $\BB$ as $\CC[\ii_2]$, as in \eqref{embed-cx}.
The equations \eqref{cx-char2} reduce to the harmonic morphism equations for maps from (an open subset of) $\RR^3$ to $\RR^2 = \CC$ and with $G =  g \in \CC = \RR^2$ and $H = h$,
\eqref{W-cx} reduces to the Weierstrass representation in
\cite{Ba-Wo-Bernstein} and \cite[(1.3.18)]{book}.
Examples \ref{ex:proj-cx}, \ref{ex:radial-cx} with
$\ve = \pm 1$ and \ref{ex:disc-cx}
reduce to the standard examples in \cite[Section 1.5]{book}. 

However, with $\ve = \pm j$, the degenerate complex-harmonic morphism of Example \ref{ex:radial-cx} does not restrict to any harmonic morphism from an open subset of $\RR^3$; indeed, all harmonic morphisms from Riemannian manifolds are non-degenerate everywhere.
 
We also have \cite{Ba-Wo-Bernstein} a Bernstein-type theorem that orthogonal projection $\RR^3 \to \RR^2$ is the only globally defined harmonic morphism from $\RR^3$ to a Riemann surface, up to isometries and postcomposition with weakly conformal maps.

The directions of fibres are parametrized by $S^2$.
The inclusion map $S^2_{\CC} \hookrightarrow \Q^2_{\CC}$ restricts to a conformal diffeomorphism of $S^2$ onto the real points $\Q^2_{\RR}$ of $\Q^2_{\CC}$, and the standard chart $\BB \to S^2_{\CC} \hookrightarrow \Q^2_{\CC}$  (Example \ref{ex:complex-quadric}(i)) restricts to the standard chart $\CC \to S^2  \stackrel{\cong}\rightarrow \Q^2_{\RR}$, exhibiting the conformal compactification of $\CC$ as $S^2$ or, equivalently, $\Q^2_{\RR}$\,. 

\bigskip 

Next, let $M^m = \RR^m_1$ be \emph{Minkowski space}, i.e.,
 $\RR^m$ endowed with the metric of signature $(1,m-1)$ given in standard coordinates $(x_1,x_2,\ldots,x_m) \in \RR^m$ by $g = - \dd x_1{}\!^2 + \dd x_2{}\!^2 + \ldots \dd x_m{}\!^2$\,.  Let $\phi:\MM^m \to \RR$ or $\CC$ be a smooth map.  Consider the following equations
\begin{equation} \left\{ \begin{array}{rrcl} \ds
\text{(a)} &\quad \ds -\frac{\pa^2\phi}{\pa x_1{}\!^2} + \sum_{i=2}^m \frac{\pa^2\phi}{\pa x_i{}\!^2}
& = & 0\,, \\
\text{(b)} & \quad \ds
-\left( \frac{\pa\phi}{\pa x_1}\right)^{\!\! 2} +
	\sum_{i=2}^m \left( \frac{\pa\phi}{\pa x_i}\right)^{\!\! 2}  & = & 0\,,
\end{array} \right. \label{null-wave}
\end{equation}
for $(x_1, \ldots, x_m) \in U$.
Then $\phi$ is harmonic if and only if it satisfies the
\emph{wave equation} (\ref{null-wave}a).  It is horizontally weakly conformal if and only if it is null in the sense that it satisfies (\ref{null-wave}b).  Hence,
$\phi$ is a harmonic morphism if and only if it satisfies both equations \eqref{null-wave}, i.e., it is a \emph{null solution of the wave equation}. 

To fit these into our theory, embed $\CC$ in $\BB$ as $\CC[\ii_2]$, and embed $\RR^3_1$ in
$\CC^3 = \CC[\ii_1]^3 \subset \BB^3$ by
$(x_1,x_2,x_3) \mapsto (x_1,x_2\ii_1, x_3\ii_1)$.  Then the equations \eqref{cx-char2} for a complex-harmonic morphism reduce to the harmonic morphism equations \eqref{null-wave}.  On setting $G = g\ii_1 $ and $H = h\ii_1 $ we obtain the Weierstrass representation obtained in
\cite[\S 2]{Ba-Wo-Rou}.

The possible directions of (non-degenerate) fibres are parametrized by the hyperbola
$H^2 = \{(x_1,x_2,x_3) \in \RR^3_1 : - x_1{}\!^2 + x_2{}\!^2 + x_3{}\!^3 = -1 \}$.
The embedding
$(x_1,x_2,x_3) \mapsto (x_1,x_2\ii_1 ,x_3\ii_1)$ maps $H^2$ 
into $S^2_{\CC}$, and thus into $\Q^2_{\CC}$ with image lying in the quadric $\bigl\{[\eta_0,\eta_1,\eta_2,\eta_3] \in \RP^3: \eta_0{}\!^2 = \eta_1{}\!^2 - \eta_2{}\!^2 - \eta_3{}\!^2\bigr\} \cong S^2$; this quadric is thus a conformal compactification of $H^2$.

As regards  Example \ref{ex:radial-cx} (complex radial projection) with
$G = q\ii_1$ and $H=0$, the solutions with $\ve = \pm 1$ restrict to radial projection from the interior of the light cone of $\MM^3$ to the hyperbola $H^2$.  On writing $j$ as $\ii_1\ii_2$ and putting the $\ii_1$ under the square root, we see that the solutions with $\ve = \pm j$ restrict to a degenerate harmonic morphism on the exterior of the light cone with fibres the tangent planes to the light cone, see \cite[Example 2.10]{Ba-Wo-Rou} for more details on these harmonic morphisms.  

The complex disc example (Example \ref{ex:disc-cx}) restricts to a globally defined surjective submersive harmonic morphism from Minkowski $3$-space $\MM^3 = \RR^3_1$ to the unit disc;  thus \emph{there is a globally defined harmonic morphism other than orthogonal projection}, in contrast to Bernstein-type theorem for the Euclidean case mentioned above.

\section{Harmonic morphisms to a Lorentz
surface}\label{sec:R31toLorentzsurface}

To discuss harmonic morphisms to a Lorentz surface, we shall use the hyperbolic numbers.  Let $\DD = \{ (x_1, x_2) \in \RR^2\}$ equipped with the usual coordinate-wise addition, but with multiplication given by
$$
(x_1, x_2)\,(y_1, y_2) = (x_1y_1+x_2y_2\,,\, x_1y_2+x_2y_1)\,.
$$
The commutative algebra $\DD$ is called the \emph{hyperbolic} (or \emph{double} or \emph{paracomplex}) \emph{numbers}. 
Write $\jj = (0,1)$; then we have $(x_1,x_2) = x_1+x_2\jj$ with $\jj^2 = 1$.  Note that $\DD$ has zero divisors, namely the numbers $a(1\pm \jj)\ (a\in \RR )$. 
By analogy with the complex numbers, we say that a $C^2$ map $\phi : U \to \DD$, $w = \phi(z)$, from an open subset of $\DD$ is \emph{H-holo\-morphic} (resp., \emph{H-anti\-holo\-morphic}) if, on setting
$z = x_1 + x_2 \jj$ and $\ov{z} = x_1 - x_2 \jj$\,, we have
$$
\frac{\pa w}{\pa \ov{z}} = 0 \quad \left( {\rm resp.,} \ \frac{\pa w}{\pa z} = 0 \right) \,;
$$
equivalently, on setting $w = u_1 + u_2 \jj$\,,  the map $\phi$ satisfies the \emph{H-Cauchy-Riemann equations}:
$$
\frac{\pa u_1}{\pa x_1} = \frac{\pa u_2}{\pa x_2} \ {\rm and} \
\frac{\pa u_1}{\pa x_2} = -\frac{\pa u_2}{\pa x_1}  \quad
\left( {\rm resp.,} \ \frac{\pa u_1}{\pa x_1}
	= \frac{\pa u_2}{\pa x_2}\ {\rm and} \ \frac{\pa u_1}{\pa x_2}
	= -\frac{\pa u_1}{\pa x_2}\right).
$$
 
By a \emph{Lorentz} surface, we mean a smooth surface equipped with a conformal equivalence class of Lorentzian metrics --- here two metrics $g, g^{\prime}$ on $N^2$ are said to be \emph{conformally equivalent} if $g^{\prime} = \mu g$ for some (smooth) function $\mu : N^2 \to \RR \setminus \{ 0\}$.  Any Lorentz surface is locally conformally equivalent to $2$-dimensional Minkowski space $\MM^2$, see, for example, \cite{book}.  
Let $\phi :U\to N^2_1$ be a $C^2$ mapping
from an open subset\/ $U$ of\/ $\RR^3_1$ to a Lorentz surface. For local considerations, we can assume that $\phi$ has values in
$\MM^2$.  Then, on identifying $\MM^2$ with the space $\DD$
of hyperbolic numbers as above, and writing $\phi =
\phi_1+\phi_2\jj$\,, the map $\phi$ is a harmonic morphism if and only if it satisfies equations \eqref{null-wave} with $m=3$, where now $\phi$ has values in $\DD$.

Now the hyperbolic numbers $\DD$ can be embedded in $\BB$ by
\begin{equation} \label{embed-hyp}
\iota_{\DD} : \DD \hookrightarrow \BB\,,  \quad \iota_{\DD} (x+y\jj ) = x + (y\ii_1) \ii_2 = x+y \jj  \quad (x,y \in \RR ) \,;
\end{equation}
this preserves all the arithmetic operations; in fact we can think of $\BB$ as the complexification $\DD \otimes_{\RR}\CC$ of $\DD$, as well as the complexification of $\CC$. Further, we have a version of Lemma \ref{lem:derivs}, as follows.

\begin{lemma}
Let $f:U \to \CC$ be real-analytic H-holo\-morphic map from an open subset of\/ $\DD$. 
Then $f$ can be extended to a bi\-com\-plex-holo\-mor\-phic function $\psi:\wt{U} \to \BB$ on an open subset $\wt{U}$ of\/ $\BB$ containing $U$; the germ of the extension at $U$ is unique.

Conversely, the restriction of any
bi\-com\-plex-holo\-mor\-phic function $\wt{U} \to \BB$ to
$U = \wt{U} \cap \DD$ is real analytic and H-holo\-morphic,
provided that $U$ is non-empty.
\end{lemma}

\begin{proof}  Write points of $U\subseteq \DD$ in the form $x+y\jj$\,; then the map $\iota_{\DD} (x+y\jj ) = q_1+q_2\ii_2$
given by $q_1 = x$ and $q_2 = y i_1$ identifies $U$ with a subset of $\BB$ which we continue to denote by $U$.  Write $f : U \to \DD$ in the form
$f(x+y\jj ) = u_1(x,y) + u_2(x,y) \jj$\,.
Extend the functions $u_i(x,y)$ by analytic continuation to holomorphic functions $u_i(q_1, q_2)$ ($i=1,2$) on an open subset $\wt{U} \supset U$ of $\CC^2 \cong \BB$ and define $\psi : \wt{U} \to \BB$ by 
$\psi (q_1+q_2\ii_2)
  = \psi_1(q_1, q_2) + \psi_2(q_1, q_2)\ii_2$
where $\psi_1 = u_1$ and $\psi_2 = u_2 \ii_1$. 
For each $i=1,2$, write $q_i = x_i + y_i\ii_1$; then since $\psi_i$ is complex analytic, on $U$ we have
\begin{eqnarray*}
\frac{\pa \psi_1}{\pa q_1} & = & \frac{\pa \psi_1}{\pa x_1} 
  =  \frac{\pa u_1}{\pa x} \quad \text{and} \\
\frac{\pa \psi_2}{\pa q_2} & = &
		- \frac{\pa \psi_2}{\pa y_2}\ii_1  
  =  - \frac{\pa}{\pa y}\ii_1(u_2 \ii_1) =  \frac{\pa u_2}{\pa y}\,.
\end{eqnarray*}
Hence, on $U$,
$$
\frac{\pa u_1}{\pa x} = \frac{\pa u_2}{\pa y}
\quad\text{if and only if}\quad
\frac{\pa \psi_1}{\pa q_1} = \frac{\pa \psi_2}{\pa q_2}\,.
$$
Similarly,
$$
\frac{\pa u_1}{\pa y} = \frac{\pa u_2}{\pa x}
\quad\text{if and only if}\quad
\frac{\pa \psi_1}{\pa q_2} = - \frac{\pa \psi_2}{\pa q_1}\,.
$$
Now, if the right-hand equations hold on $U$ then, by analytic continuation,  they hold on $\wt{U}$ proving the first part of the lemma; the converse is similar.  \end{proof}

To recover the formulae of \cite[\S 3]{Ba-Wo-Rou} for harmonic morphisms from
$\MM^3 = \RR^3_1$ to $\MM^2 = \DD$, embed $\RR^3_1$ in $\CC^3 = \CC^3[\ii_1] \subset \BB^3$  by $(x_1,x_2,x_3) \mapsto (x_3, x_1\ii_1, -x_2)$ --- this is a different embedding to that used in in \S \ref{sec:real}.  Nondegenerate fibres are now spacelike lines whose directions are parametrized by the pseudosphere
$S^2_1 = \{(x_1,x_2,x_3) \in \RR^3_1: - x_1{}\!^2 + x_2{}\!^2 + x_3{}\!^3 = 1$.  This is mapped into $S^2_{\CC}$, and thus into $\Q^2_{\CC}$, with image in the quadric
$\bigl\{[\zeta_0,\zeta_1,\zeta_2,\zeta_3] \in \RP^3: \zeta_0{}\!^2
= \zeta_2{}\!^2 + \zeta_3{}\!^2 - \zeta_1{}\!^2\bigr\} \cong S^1 \times S^1$.
This quadric is the standard conformal compactification of $S^2_1$ and of $\MM^2$, see \cite[Example 14.1.2]{book} for more details.  Then set $G = g\ii_1 $ and $H = h\ii_1 $.

As regards  Example \ref{ex:radial-cx} (complex radial projection) with
$G = q\ii_1$, the solutions with $\ve = \pm 1$ restrict to radial projection from the \emph{exterior} of the light cone of $\MM^3$ to the pseudosphere $S^2_1$. The solutions with $\ve = \pm j$ restrict to a degenerate harmonic morphism, again defined on the exterior of the light cone, with fibres the tangent planes to the light cone; see \cite[Example 3.5]{Ba-Wo-Rou} for more details on these harmonic morphisms. 

On setting $t = \ii_1$, the complex disc example (Example \ref{ex:disc-cx}) restricts to a harmonic morphism from an open subset of $\MM^3$, see \cite[Example 3.6]{Ba-Wo-Rou} for a description.

\end{document}